\documentclass[final,1p,times]{elsarticle}

\usepackage{amssymb}

\newcommand{\defeq}{\mathrel{\mathop:}=}

\usepackage{geometry}
\usepackage{hyperref}
\usepackage{svg}
\usepackage{amsmath}
\usepackage{amsfonts}
\usepackage{subcaption}
\usepackage[ruled,lined]{algorithm2e}
\usepackage{multirow}
\usepackage{xr}
\usepackage{listings}
\usepackage{hyperref}
\usepackage[numbers]{natbib}

\journal{Applied Mathematics and Computation}

\begin{document}

\begin{frontmatter}

\title{A C++ implementation of the discrete adjoint sensitivity analysis method for explicit adaptive Runge-Kutta methods enabled by automatic adjoint differentiation and SIMD vectorization}

\author[affil1]{Rui Martins}
\ead{rui.carlos.andrade.martins@gmail.com}
\author[affil1]{Evgeny Lakshtanov}
\ead{lakshtanov@ua.pt}
\affiliation[affil1]{organization={Department of Mathematics,University of Aveiro},
addressline={Campus Universitário de Santiago},
postcode={3810-193},
city={Aveiro},
country={Portugal}}

\begin{abstract}
A C++ library for sensitivity analysis of optimisation problems involving ordinary differential equations (ODEs) enabled by automatic differentiation (AD) and SIMD (Single Instruction, Multiple data) vectorization is presented. The discrete adjoint sensitivity analysis method is implemented for adaptive explicit Runge-Kutta (ERK) methods. Automatic adjoint differentiation (AAD) is employed for efficient evaluations of products of vectors and the Jacobian matrix of the right hand side of the ODE system. This approach avoids the low-level drawbacks of the black box approach of employing AAD on the entire ODE solver and opens the possibility to leverage  parallelization. SIMD vectorization is employed to compute the vector-Jacobian products concurrently. We study the performance of other methods and implementations of sensitivity analysis and we find that our algorithm presents a small advantage compared to equivalent existing software.
\end{abstract}


% 34-04 Software, source code, etc. for problems pertaining to ordinary differential equations
%65L06 Multistep, Runge-Kutta and extrapolation methods for ordinary differential equations
%65K10 Numerical optimization and variational techniques
%90C31 Sensitivity, stability, parametric optimization

\end{frontmatter}

\section{Introduction} \label{sec::introduction}
Sensitivity analysis focuses on quantitatively evaluating how the outputs of a mathematical model respond to variations in its inputs. In local sensitivity analysis we measure how an individual small change in the vicinity of one of the model's parameters affects its output. In the context of differential equations, sensitivities are effectively the derivatives of an objective functional which depends on the model's trajectory with respect to either the parameters or the initial conditions of the problem. Applications of sensitivity analysis are found in parameter estimation \cite{navon1998practical,zhang2014parameter}, chemical kinematics \cite{KPP2003}, systems biology \cite{sommer2017numerical}, ocean and atmosphere dynamics \cite{sandu2005adjoint}, dynamic optimisation \cite{ozyurt2005large}, data assimilation \cite{le2002second}, optimal control \cite{griesse2003parametric}, and recently neural networks associated with ODEs \cite{chen2018neural,grathwohl2018ffjord,poli2019graph}.

In this context, employing conventional methods for numerically evaluating derivatives is either inefficient or inaccurate. Numerical differentiation is prone to truncation errors and floating point cancellation errors \cite{burden2011numerical}, symbolic differentiation is memory intensive, sluggish and cannot handle some control-flow mechanisms \cite{margossian2019review} and a black-box approach using automatic differentiation produces non-optimal code and suffers from low-level drawbacks, such as memory allocation, management of pointers, I/O and integration of parallelization \cite{petsc2022}. 

The most sophisticated approaches to sensitivity analysis of ODEs are classified in two distinct classes and involve a mixture of hand-written code with one of the derivative evaluation techniques mentioned in the previous paragraph: the continuous approach (differentiate-then-discretize), where the original system of ODEs is augmented with the sensitivities' dynamics, derived by hand, and only afterwards is solved numerically; the discrete approach (discretize-then-differentiate), where we discretize the governing equations with an appropriate numerical method and then differentiate the numerical scheme itself to obtain a scheme for computing the sensitivities. Additionally, one can distinguish between forward and adjoint mode differentiation. In forward mode differentiation, we compute derivatives of intermediate variables with respect to a fixed independent variable (effectively one of the parameters or initial condition), obtaining the sensitivities of the model outputs with respect to the independent variable in a single algorithm run, while in adjoint mode differentiation we compute derivatives of the output with respect to the intermediate variables, obtaining the gradient of one of the model's outputs with respect to all the independent variables. These approaches may differ widely in implementation, each with its own benefits and counterparts. In state of art sensitivity analysis software, these methods leverage AD for evaluating vector-Jacobian or Jacobian-vector products efficiently, although many also support finite-differences coloring schemes for the same purpose.

Of the several tools that have been developed to implement sensitivity analysis for ODEs efficiently and automatically, we now make an effort to enumerate the most significant in the literature. The FORTRAN package ODESSA \cite{odessa1988}, which is for continuous forward sensitivity analysis of ODEs using the LSODE (Livermore Solver for Ordinary Differential Equations) solver. The C package CVODES, within the SUNDIALS software suite \cite{cvodes2005}, is an ODE solver implementing continuous forward and adjoint sensitivity analysis for adaptive backward differentiation formulas and Adams–Moulton methods of variable order. The Kinectic Pre-Processor (KPP) \cite{KPP2003} generates FORTRAN and MATLAB code that implements continuous forward and adjoint sensitivity analysis and discrete adjoint sensitivity analysis for backward differentiation formulas and Rosenbrock methods, and is specifically aimed at chemical kinetic systems. The DASPK library in \cite{daspk1999} supports continuous sensitivity analysis for differential-algebraic equations solved by adaptive backward differentiation formulas with parallel support. In \cite{chen2018neural} the authors propose an implementation of the continuous adjoint method for training continuous depth neural networks, but this approach is sometimes unstable and prone to large numerical errors since the integration of the adjoint system does not take into account the forward trajectory \cite{juliaarticle} (it functions well enough in their case given that for gradient descent, a rough estimate of the gradient is good enough). In the FORTRAN library DENSERKS \cite{DENSERKS2009}, the continuous adjoint sensitivity method is implemented for ERK methods, leveraging the built in dense output mechanism of ERK implementations to accurately reconstruct the Trajectory at any time instance. This approach offers flexibility compared to the discrete adjoint method, since the time steps taken during the adjoint ODE solve may differ from those taken during the forward run, and may be adapted according to the user's accuracy needs, but at extra computational cost. The NIXE library \cite{nixelotz2015higher} is a discrete adjoint ODE integration framework written in C++ combined with operator overloading automatic differentiation (dco/c++) for linearly-implicit Euler discretizations. The FORTRAN library FATODE \cite{fatode2014} is a general purpose package that provides discrete forward and adjoint sensitivity analysis for explicit Runge-Kutta, fully implicit Runge-Kutta, singly-diagonal Runge-Kutta methods and Rosenbrock methods. The dolfin-adjoint (or pyadjoint) package \cite{mitusch2019dolfin} employs operator-overloading automatic differentiation to generate code for forward and adjoint sensitivity analysis at runtime, aimed at finite element models in the FEniCS software environment, supporting parallelization. Among recently developed sensitivity analysis software stands the DifferentialEquations.jl \cite{differentialequationsjl2017} and SciMLSensitivity.jl \cite{rackauckas2020universal}, general purpose packages for solving differential equations and for sensitivity analysis in the Julia programming language, with support for continuous first-order forward and adjoint sensitivity analysis for ODEs and for discrete sensitivity analysis via using AD on the entire solver. The 2024 release of Matlab introduced sensitivity analysis for ordinary differential equations using the CVODES software suite. Finally, the TSAdjoint component in the Portable Extensible Toolkit for Scientific Computation (PETSC) \cite{petsc2022} is a general purpose HPC-friendly library for first and second-order discrete adjoint sensitivity analysis of several time integration methods for ODEs and differential-algebraic equations.

This article presents an efficient C++ implementation of sensitivity analysis for optimisation problems involving systems of first-order ODEs, implementing the discrete adjoint sensitivity method for adaptive ERK methods available in the C++ \texttt{boost::odeint} library. Automatic adjoint differentiation is used to evaluate products between vectors and the Jacobian of the right-hand-side of the system of ODEs using the C++ package of \cite{Matlogica} named AADC. The AADC library also allows for SIMD vectorization of the vector-Jacobian products, accelerating the computation of the sensitivities for multiple objective functions.

To the best of the authors' knowledge, FATODE and TSAdjoint are the only libraries implementing the discrete adjoint sensitivity analysis method for adaptive ERK methods. In contrast, DENSERKS implements the continuous adjoint sensitivity analysis method for adaptive ERK methods. In this article we propose a couple of changes on these implementations that we proceed to enumerate. The first change is the vectorization of the vector-Jacobian products via leveraging the SIMD functionality of the AADC framework, allowing the computation of sensitivities of multiple objective functions in parallel, which none of the above mentioned implementations currently figures. The second is the implementation of the computation of vector-Jacobian products efficiently using AAD, without the need to explicitly evaluate the Jacobian matrices, something which has also been addressed by DENSERKS and TSAdjoint \footnote{To the best of the authors' knowledge, the version of FATODE being considered evaluates the Jacobian matrix first and only then computes the vector-Jacobian products}. FATODE saves the forward trajectory of the system of ODEs in a tape, without implementing a check-pointing scheme. Our C++ library uses a simple check-pointing scheme to alleviate the memory costs, DENSERKS uses a two level check-pointing scheme, and TSAdjoint uses a sophisticated online optimal check-pointing scheme that accounts for the memory limitations of specific hardware. We remark that in TSAdjoint, the adjoint models are implemented based on the parallel infrastructure in PETSc, whereas this kind of parallelization is not leveraged to accelerate the adjoint solver in our implementation. However, the parallelization of sensitivity analysis is not automatic in PETSc (the user needs to implement it correctly), whereas our SIMD approach using Advanced Vector Extensions (AVX) to evaluate the sensitivities of several objective functions simultaneously requires the user no extra work.

In section \ref{sec:modelproblem} we establish the mathematical foundation of the proposed algorithm by presenting the model system of ODEs for which we will perform sensitivity analysis and providing a description of adaptive $s$-stage ERK methods. We proceed by providing a general description of some methods of sensitivity analysis of ODEs in section \ref{sec::adjoints} and deriving the discrete adjoint sensitivity analysis scheme for general adaptive ERK schemes. In section \ref{sec:results} we validate our algorithm by computing the sensitivities of a two-dimensional heat equation initial value problem, where the sensitivities are known analytically. Additionally, we perform benchmarks by comparing the performance of our algorithm to other methods and implementations of sensitivity analysis on the generalised Lotka-Volterra equations. Section \ref{sec:conclusion} presents the conclusion of this work.

\section{Model Problem} \label{sec:modelproblem}
Consider the vector $\mathbf{\alpha} = (\alpha_1,\dots,\alpha_P)^T$ of $P\in\mathbb{N}$ real valued parameters. In this article, we will be considering an initial value problem for a system of $N\in\mathbb{N}$ parameter dependent ODEs in the interval $[t_0,t_f]$. We call $\mathbf{u}(t) = \left(u_1(t),\dots,u_N(t)\right)^T$ the state vector, representing the solution at time $t$. Our model problem is therefore given by:

\begin{align}
    \frac{d \mathbf{u}(t)}{d t} &= F(\mathbf{u}(t),\alpha,t), \qquad
    \mathbf{u}(t_0) = \mathbf{u}^0
    \label{eq::system_of_odes}
\end{align}
Here, $F\colon\mathbb{R}^N\times\mathbb{R}^P\times\mathbb{R}\to\mathbb{R}$ is a real-valued $C^1$ function, $\{u^0_i\}_{i=1}^N$ are real numbers representing the initial conditions and $t$ is the independent variable (time). Consider, additionally, $M\in\mathbb{N}$ objective/cost functionals $\psi_i$ of the form

\begin{equation}
    \psi_i\left(\mathbf{u}(t_f),\mathbf{u}(t_0),\alpha\right) = E_i(\mathbf{u}(t_0),\mathbf{u}(t_f),\alpha)+\int_{t_0}^{t_f} R_i(\mathbf{u}(t),\alpha,t).
    \label{eq::objective_functional}
\end{equation}
The first term $E_i$ is referred to as an end-point cost, while the second integral is referred to as a trajectory cost. Here, $E_i\colon\mathbb{R}^N\times\mathbb{R}^N\times\mathbb{R}^P\to\mathbb{R}$ and $R_i\colon \mathbb{R}^N\times\mathbb{R}^P\times\mathbb{R}\to\mathbb{R}$ are real-valued scalar $C^1$ functions. The aim of sensitivity analysis for ODEs is to efficiently compute the following total derivatives:

\begin{equation}
    \frac{d \psi_i}{d u^0_j}, \qquad \frac{d \psi_i}{d \alpha_k},
\end{equation}
for which we will employ the discrete adjoint sensitivity analysis method enhanced with automatic differentiation and vectorization.

\subsection{The adaptive explicit Runge-Kutta Method\label{section::runge_kutta_discretization}} 
Consider a discretization $(t^n)_{n=0}^T$ of the time interval $[t_0,t_f]$, with $T\in\mathbb{N}$, and define the time steps $\Delta t^n \defeq t^{n+1}-t^n$. The approximation of the solution at time $t = t^n$ will be $\mathbf{u}^n \approx \mathbf{u}(t^n)$. An adaptive explicit $s$-stage Runge-Kutta method to compute the approximate solution $\mathbf{u}^T\approx\mathbf{u}(t_f)$ of the initial-value problem in \ref{eq::system_of_odes} is given by the following iterative formula, which can be found on \cite{lambert1991numerical}:

\begin{equation}
	\mathbf{u}^{n+1} = \mathbf{u}^n+\Delta t^n\sum_{m=1}^s  b_m \mathbf{k}^m
    \label{eq::rungekutta1}
\end{equation}
where the so-called $\mathbf{k}$-vectors are given by,

\begin{equation}
	\mathbf{k}^m = F\big(\mathbf{u}^{m,n},\alpha,t^{m,n}\big)
    \label{eq::rungekutta2},
\end{equation}
with intermediate Runge-Kutta state,

\begin{equation}
    \mathbf{u}^{m,n}=\mathbf{u}^n+\Delta t^n\sum_{j=1}^{m-1} a_{mj}\mathbf{k}^j,
\end{equation}
and intermediate time $t^{m,n}$,

\begin{equation}
    t^{m,n} = t^n +c_m\Delta t^n.
\end{equation}
To specify a particular method, we need to pick the number of stages $s$, the real-valued explicit Runge-Kutta matrix $(a_{i j})_{j<i}^s$, the weights $(b_i)_{i=1}^s$ and the nodes $(c_i)_{i=1}^s$, which are usually compiled in a Butcher Tableau \cite{butcher2016numerical}. 

One can implement adaptive time-stepping at negligible computational cost by computing an estimate of the local truncation error for a single Runge-Kutta step. This estimate is given by comparing the numerical solution $\mathbf{u}^{n+1}$ obtained from the Runge-Kutta method of order $O_S$ with another solution $\mathbf{u}_*^{n+1}$ obtained from a method of lower order $O_E$, in all things equal to its higher order counterpart except in the choice of weights. The local error can thus be estimated as follows:

\begin{equation}
	\epsilon^n_i \approx u_{*i}^{n+1}-u_i^{n+1}
	= \Delta t^n \sum_{m=1}^s (b_m-b_m^*)k^m_i,
\end{equation}
where $(b_m^*)_{m=1}^s$ are the weights of the lower order Runge-Kutta method. Originally proposed by Fehlberg \cite{fehlberg1969low}, this class of methods are known as embedded Runge-Kutta methods. 

This estimate can be used to assert if the current time-step $\Delta t^{n}$ keeps the truncation error in control relative to an absolute $\epsilon_a$ and relative $\epsilon_r$ user-defined bounds. In the C++ library \texttt{odeint} \cite{odeint}, the time-step control is done as follows. Let,

\begin{equation}
    v_i =\frac{|\epsilon^n_i|}{\epsilon_a+\epsilon_r(|u_i^n|+|\frac{du^n_i}{dt}|)}
    \qquad \text{and}\qquad
    v = \left\lVert\mathbf{v}\right\rVert,
\end{equation}
where $||\cdot||$ is a vector norm (usually, the Euclidean norm or inf-norm) and 

\begin{equation}
    \frac{d\mathbf{u}^n}{dt} = \Delta t^n\sum_{m=1}^s b_m \mathbf{k}^m. 
\end{equation} 
Then, the new time-step $\Delta t^{'}$ is given by:

\begin{equation}
    \begin{cases}
    \Delta t^\prime =  \Delta t^n \max(\frac{0.9}{v^{O_E-1}},0.2) \qquad \text{if } v>1\\
    \Delta t^\prime =  \Delta t^n \min(\frac{0.9}{v^{O_S}},5) \qquad \text{if } v<0.5\\
    \Delta t^\prime =  \Delta t^n \ \qquad \text{otherwise}
    \end{cases}
\end{equation}
If $v>1$, the original time step $\Delta t^n$ was to big. When this happens, $\Delta t^\prime$ is used to re-compute the solution at iteration $n+1$. This process is repeated until the error falls within the absolute and relative tolerances provided by the user, and in the end we take $\Delta t^{n+1}=\Delta t^\prime$. In contrast, if $v<1$, we use $\Delta t^\prime$ directly as the time-step for the next iteration, $\Delta t^{n+1}=\Delta t^{\prime}$.

\section{Description of some sensitivity analysis algorithms\label{sec::adjoints}}
In this section we provide a general description of some methods of sensitivity analysis for ODEs and describe the theory behind the discrete adjoint sensitivity analysis method in particular, which was implemented in this work.

\subsection{Numerical differentiation or finite-differences method} 
A straightforward approach to sensitivity analysis is numerical differentiation (ND) where we approximate the sensitivities via an expression based on the Taylor series expansion of the objective function, such as the forward finite differences:

\begin{equation}
    \frac{d \psi_i}{d\alpha _k} \approx \frac{\psi_i(\mathbf{u}_{\alpha+\Delta\alpha_k}(t),\alpha+\Delta\alpha_k)-\psi_i(\mathbf{u}_{\alpha}(t),\alpha)}{\Delta\alpha_k}
    \label{eq::numerical_differentiation_expression}
\end{equation}
where $\mathbf{u}_{\alpha}(t)$ is the numerical trajectory of the ODE system with initial conditions $\mathbf{u}^0$ and parameters $\alpha$ and ${\Delta\alpha_k}$ is a vector of zeroes except for the $k$th entry, where it is equal to a small value $\Delta\alpha$ called the step or bump (a similar expression to \ref{eq::numerical_differentiation_expression} can be used for the derivatives with respect to the initial conditions). This method requires evaluating the solution of the system of ODEs twice for each parameter sensitivity, making it of time complexity of $\mathcal{O}(N\times P)$ when computing the sensitivities of $\psi_i$ for a single index $i$. Additionally, the method is prone to truncation errors and finite-precision cancellation errors.

\subsection{Automatic Differentiation}
Automatic differentiation can be directly employed for sensitivity analysis of simple models described by a computer program, requiring little technical knowledge and few code changes from the user. However, employing a black-box approach of AD on an ODE solver produces non-optimal code that suffers from low-level drawbacks, and the memory requirements might be very large \cite{fatode2014,petsc2022}. Instead, in state-of-the-art sensitivity analysis software, automatic differentiation is used in conjunction with regular forward or adjoint differentiation formulas for the numerical methods considered, often used to evaluate products between vectors and the Jacobian of the right-hand-side, either it be in its continuous or discrete setting.

\subsubsection{Continuous forward sensitivity analysis method}
In the continuous forward-mode sensitivity analysis method (CFSA), the original system of ODEs is augmented with the following dynamic equations (Einstein's notation is adopted hereafter for the sub-indices unless specified):

\begin{align}
    \dot{\gamma}_{ij}(t) &= J^x_{il}(\mathbf{u}(t),\alpha,t)\gamma_{lj}(t) \label{eq::CFSAinitialconditions}
    \\
    \dot{\beta}_{ik}(t) &= J^x_{ij}(\mathbf{u}(t),\alpha,t)\beta_{jk}(t)+J^\alpha_{ik}(\mathbf{u}(t),\alpha,t)\label{eq::CFSAparameters}
    \\ 
    \gamma_{ij}(t_0) &= \delta_{ij},\qquad \beta_{ik}(t_0) = 0
    \label{eq::CFSA_ODE}
\end{align}
where $\beta_{ik}(t)$ are the derivatives of $u(t)$ with respect to the parameters and $\gamma_{ij}(t)$ the derivatives of the $u(t)$ with respect to the initial conditions:

\begin{equation}
        \gamma_{ij}(t) = \frac{d u_i(t)}{d u^0_j}, \qquad \beta_{ik}(t) = \frac{d u_i(t)}{d\alpha_k},
\end{equation}
and $J^x_{ij}$ and $J^\alpha_{ik}$ are the Jacobian matrices of the right-hand-side $F$ with respect to the state vector and parameters,
\begin{align}
    J^x_{ij}(\mathbf{u},\alpha,t) = \frac{\partial F_i(\mathbf{u},\alpha,t)}{\partial u_j},
    \\
    J^\alpha_{ik}(\mathbf{u},\alpha,t) = \frac{\partial F_i(\mathbf{u},\alpha,t)}{\partial \alpha_k},
\end{align}
When employing this method the objective is usually to compute the sensitivities of the solution of the original ODE system at time $t_f$ with respect to the initial conditions and parameters, which will be given by:

\begin{equation}
    \frac{d u_i(t_f)}{du^0_j} = \gamma_{ij}(t_f) \qquad \frac{d u_i(t_f)}{d\alpha_k} = \beta_{ik}(t_f) 
\end{equation}
If one is interested in only the derivatives with respect to the initial conditions, one would augment the original system only with \ref{eq::CFSAinitialconditions}. Alternatively, if one is interested only in the derivatives with respect to the parameters, one would augment the original system with \ref{eq::CFSAparameters}. In an efficient approach, the Jacobian matrix is not explicitly computed, and instead, for each parameter $k$, the Jacobian-vector products $ J^x_{ij}(\mathbf{u}(t),\alpha,t) v_{jk}$ can either be manually derived, evaluated via finite-differences or by \textbf{forward-mode} automatic differentiation. The time complexity for finding all the parameter sensitivities is $\mathcal{O}(N\times P)$, in which case CFSA is not practical for a large number of parameters ($P \gg N$).

\subsubsection{Continuous adjoint sensitivity analysis method} 
Another mode of continuous sensitivity analysis, called the continuous adjoint sensitivity analysis (CASA) method, is better suited to find the sensitivities for optimisation problems involving a large number of parameters. Usually, the aim is to compute the total derivatives of objective functionals in \ref{eq::objective_functional} with respect to the initial conditions and parameters, but in this section we focus only on the parameter sensitivities for the sake of brevity. We will treat the end-point and trajectory terms separately, focusing on the trajectory term first:

\begin{equation}
    \psi_i^{(2)}[\mathbf{u}(t),\alpha] = \int_{t_0}^{t_f} R_i(\mathbf{u}(t),\alpha,t)
\end{equation}
Following the approach with Lagrange multipliers described in \cite{DENSERKS2009}, we define the following set of adjoint variables described by ODEs with final value  condition:

\begin{equation}
    \dot{\lambda}_{ij} = -\frac{\partial R_i}{\partial u_j}-\lambda_{il}\frac{\partial F_l}{\partial u_j}, \qquad \lambda_{ij}(t_f) = 0,
    \label{eq::adjoint_continuous}
\end{equation}
such that the total derivative of the objective with respect to the parameters can be expressed as:

\begin{equation}
    \frac{d \psi_i^{(2)}}{d\alpha_k} = \int_{t_0}^{t_f} \big(\frac{\partial R_i}{\partial \alpha_k}+\lambda_{ij}\frac{\partial F_j}{\partial \alpha_k} \big) dt +\lambda_{ij}(t_0)\frac{\partial u_j(t_0)}{\partial \alpha_k}.
\end{equation}
To compute this integral numerically, we define quadrature variables $w_{ik}$:
\begin{equation}
    \dot{w}_{ik} = -\frac{\partial R_i}{\partial \alpha_k}-\lambda_{ij}\frac{\partial F_j}{\partial \alpha_k}, \qquad w_{ik}(t_f)=0,
    \label{eq::quadrature_w}
\end{equation}
and thus the sensitivity is written as,

\begin{equation}
     \frac{d \psi_i^{(2)}}{d\alpha_k} = w_{ik}(t_0) +\lambda_{ij}(t_0)\frac{\partial u_j(t_0)}{\partial \alpha_k}.
\end{equation}
In practice, we need to solve the original system of ODEs once. Afterwards, the coupled system of $N+P$ ODEs, defined by equations \ref{eq::adjoint_continuous} and \ref{eq::quadrature_w}, must be solved backwards in time in the interval $[t_0,t_f]$ for each entry $i$ of the objective functional. Since these equations have terms that depend on the forward trajectory of the system, namely in the evaluation of the Jacobian matrices, we require the knowledge of the forward trajectory of the system. If the numerical method to solve the coupled system defined by \ref{eq::adjoint_continuous} and \ref{eq::quadrature_w} uses the same time-steps as the numerical method used to solve the original system of ODEs, then the forward trajectory must be cached, a memory intensive process. If we allow the time-steps for the solver of the coupled system to be adaptive, then the original system of ODEs must either be solved backwards in time in $[t_0,t_f]$ alongside the coupled system for the adjoints or the forward trajectory must be reconstructed as accurately as possible. For example, in \cite{DENSERKS2009}, the authors use a check-pointing scheme and the dense-output mechanism of ERK methods to achieve this.

The sensitivities for the end-point cost,

\begin{equation}
    \frac{d \psi_i^{(1)}}{d\alpha_k} = \frac{d E_i}{d\alpha_k},
\end{equation}
are easier to evaluate. By defining the following adjoint variables via the following ODEs with final value condition:

\begin{equation}
    \dot{\lambda}_{ij} = -\lambda_{il}\frac{\partial F_l}{\partial u_j}, \qquad \lambda_{ij}(t_f) = \frac{\partial E_i}{\partial u_j(t_f)},
    \label{eq::lambda2}
\end{equation}
the total derivative of the end-point objective with respect to the parameters is given by:

\begin{equation}
\frac{d \psi_i^{(1)}}{d\alpha_k} = \frac{\partial E_i}{\partial \alpha_k} + \frac{\partial E_i}{\partial u_j(t_0)}\frac{\partial u_j(t_0)}{d\alpha_k}+\lambda_{ij}(t_0)\frac{\partial u_j(t_0)}{\partial \alpha_k} + \int_{t_0}^{t_f} \lambda_{ij}\frac{\partial F_j}{\partial \alpha_k}dt 
\end{equation}
which requires solving \ref{eq::lambda2} backwards in time in the interval $[t_0,t_f]$ and evaluating the integral via a quadrature rule. 

Since the integrals above must be evaluated via a quadrature rule, this introduces some numerical error. When the forward trajectory is reconstructed inaccurately, we risk the method becoming unstable and becoming prone to large numerical errors, as happens in \cite{chen2018neural}. Finally, the vector-Jacobian products $v_{il}(t)J^\alpha_{lk}(\mathbf{u}(t),\alpha,t)$ can be either manually derived, evaluated with finite-differences or by using \textbf{reverse-mode} automatic differentiation. The time complexity of this method is $\mathcal{O}(N+P)$ for each index $i$ of the objective function, scaling better than its forward counterpart for problems with a large number of parameters ($P\gg N$).

\subsection{Discrete sensitivity analysis methods}
An alternative approach to sensitivity analysis is to directly differentiate the steps of a numerical solver in order to obtain update equations for either the derivatives of the intermediate states with respect to the parameters (in forward mode differentiation) or for the derivatives of the objective function with respect to the parameters (in adjoint mode differentiation). Even though this might be a laborious process, requiring the derivation of the solver dependant update equations by hand, in principle the method executes faster than its continuous alternatives and outputs the exact derivative of the solver program up to machine precision. Essentially, this is equivalent to applying AD to the entire solver program, but avoids the low-level problems inherent to this approach, results in better optimised code and allows for the integration of parallelization \cite{petsc2022}.

Firstly, we need to consider a discretization of the objective function to be able to compute the trajectory term numerically with the discrete solution. For that, still in the continuous setting, we consider quadrature variables $q_i(t)$ such that $\dot{q}_i(t) =  R_i(\mathbf{u}(t),\alpha,t)$ and $q_i(t_0)=0$. We can express the objective function as:

\begin{equation}
    \psi_i(\mathbf{u}(t_0),\mathbf{u}(t_f),\alpha) = E_i(\mathbf{u}(t_0),\mathbf{u}(t_f),\alpha)+q_i(t_f),
\end{equation}
as done in \cite{fatode2014}. For each entry $i$ of the objective function, we now handle the coupled augmented system of $N+1$ ODEs:

\begin{align}
\begin{bmatrix}
    \dot{\mathbf{u}}  
    \\
    \dot{\mathbf{q}}
\end{bmatrix}
=
\begin{bmatrix}
    F(\mathbf{u}(t),\alpha,t)
    \\
    R(\mathbf{u}(t),\alpha,t)
\end{bmatrix}
\end{align}
where $F$ is a $N \times 1$ vector and $R$ is a $M \times 1$ vector. In the following discussion, we will focus on computing the sensitivities for the end-point cost and the trajectory cost separately.

\subsubsection{Sensitivities of the end-point cost}
Here, we consider an objective function with no trajectory cost:

\begin{equation}
    \psi_i^{(1)}(\mathbf{u}^0,\mathbf{u}^T,\alpha) = E_i(\mathbf{u}^0,\mathbf{u}^T,\alpha)
\end{equation}
To proceed with deriving update equations for the derivatives, it is convenient that we express the solver scheme in \ref{eq::rungekutta1}, \ref{eq::rungekutta2} as a sequence of vector operations. So, for each time step $n$, we define $s+2$ intermediate state vectors $\mathbf{w}^{m,n}$, corresponding to the Runge-Kutta hidden states $\mathbf{u}^n$ and Runge-Kutta $\mathbf{k}$-vectors:

\begin{align*}
    \mathbf{w}^{0,0} &= \mathbf{u}^0, \qquad \text{if $n=0$},
    \\
    \mathbf{w}^{m,n} &= \mathbf{k}^m = F\big(\mathbf{u}^{mn},\alpha,t^{mn}\big) \qquad \text{if $m=1,\dots,s$}.
     \\
     \mathbf{w}^{s+1,n} &= \mathbf{w}^{0,n} + \Delta t^n \sum_{l=1}^s b_l \mathbf{w}^{l,n},
     \\
     \mathbf{w}^{0,n+1} &= \mathbf{w}^{s+1,n}, \qquad \text{if $n>0$ and $n<T-1$},
\end{align*}
and to each entry of these vectors $w_j^{m,n}$ we call intermediate variable. Finally, the last step is to compute the objective function:

\begin{equation}
    w^{0,T}_i = E_i(\mathbf{w}^{0,0},\mathbf{w}^{s+1,T-1},\alpha)
\end{equation}
Computing the solution with the Runge-Kutta method amounts to computing the intermediate states sequentially, where the value of the first intermediate state $\mathbf{w}^{0,0}$ corresponds to the algorithm inputs state vector $\mathbf{u^0}$, the value of the intermediate state $\mathbf{w}^{s+1,T-1}$ corresponds to the output state vector $\mathbf{u}^T$ and the value of the intermediate state $\mathbf{w}^{0,T}$ corresponds to the evaluation of the objective functions. The execution of this algorithm is referred to as the \textbf{forward pass}, and the procedure can be cast into a directed computational graph. Each node of this graph holds a value and an operation. Nodes can have parents, connected to them via an incoming edge, and children, connected to them via an outgoing edge. To start the forward pass, we load the values of the orphan nodes $w^{0,0}_i$ with the algorithm inputs $u_i^0$. To compute the value of a child node we perform the operation of that child on the values of the parent nodes. We repeat this procedure, progressing in topological order, until we reach the childless nodes $w^{0,T}_i$, which correspond to the algorithm outputs (he evaluation of the objective function).

We will skip talking about the discrete forward sensitivity method since its lengthy description would fall out of scope of this article and thus proceed directly to describe the discrete adjoint sensitivity method. Consider that to each intermediate variable $w_j^{m,n}$ and to each output vector entry $w^{0,T}_i$ corresponds an adjoint variable $\bar{w}_{ij}^{m,n}$. We call these adjoint variables with respect to the intermediate state $(m,n)$. Furthermore, to each parameter $\alpha_k$ and to each output vector entry $w^{0,T}_i$ corresponds another adjoint variable $\bar{\alpha}_{ik}$, which we call adjoint variable with respect to the $k$th parameter. We define the adjoint variables as the following total derivatives:

\begin{equation}
\bar{w}^{m,n}_{ij} \defeq \frac{d w^{0,T}_i}{d w^{m,n}_j},
\qquad
\bar{\alpha}_{ik}\defeq \frac{d w^{0,T}_i}{d\alpha_k}.
\end{equation}
To start computing the adjoint variables with respect to every node, we require the knowledge of the so-called seed, which is the adjoint variable corresponding to the output of the objective function,

\begin{equation}
    \bar{w}^{0,T}_{ij} = \delta_{ij}.
\end{equation}
To compute $\bar{w}^{s+1,T-1}_{ij}$ immediately after we require the partial derivatives of the objective function with respect to the final state of the system of ODEs, $\frac{\partial \psi_i}{\partial u^T_j}$. Furthermore, the partial derivative of the objective function with respect to the parameters $\frac{\partial \psi_i}{\partial\alpha_k}$ is required to initialise the adjoint $\bar{\alpha}_{ik}$.

Notice that if we set $\bar{w}^{s+1,T-1}_{ij}= \delta_{ij}$, then the adjoint variables $\bar{w}^{m,n}_{ij}$ correspond to the derivatives of the solution of the ODE system with respect to the intermediate states, and after a reverse pass for each output entry $i$ we would obtain the sensitivity matrix with respect to the initial conditions and with respect to the parameters.

\begin{equation}
	\bar{w}^{0,0}_{ij} = \frac{du^{T}_i}{d u^{0}_j} 
    \qquad 
    \bar{\alpha}_{ik}=\frac{du^{T}_i}{d\alpha_k},
\end{equation}
but in general, 

\begin{equation}
	\bar{w}^{0,0}_{ij} =  \frac{d \psi_i^{(1)}}{d u^{0}_j} 
    \qquad 
    \bar{\alpha}_{ik}= \frac{d \psi_i^{(1)}}{d\alpha_k},
\end{equation}
The partial derivatives of $\psi_i$ can be provided analytically, evaluated via finite-differences or via AD.

\begin{figure*}
\begin{subfigure}{.33\textwidth}
  \centering
  \includegraphics[width=.8\linewidth]{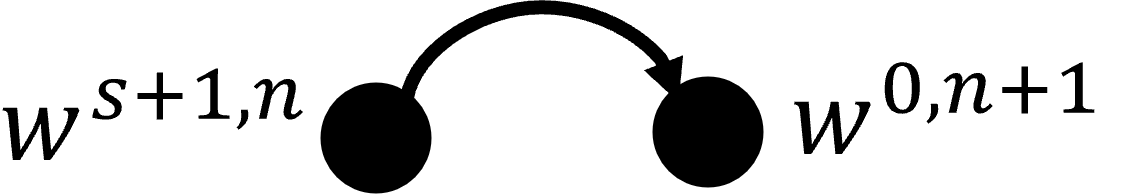}
  \caption{Corresponds to equation \ref{eq::adj1}}
  \label{fig:sfig1}
\end{subfigure}%
\begin{subfigure}{.33\textwidth}
  \centering
  \includegraphics[width=.8\linewidth]{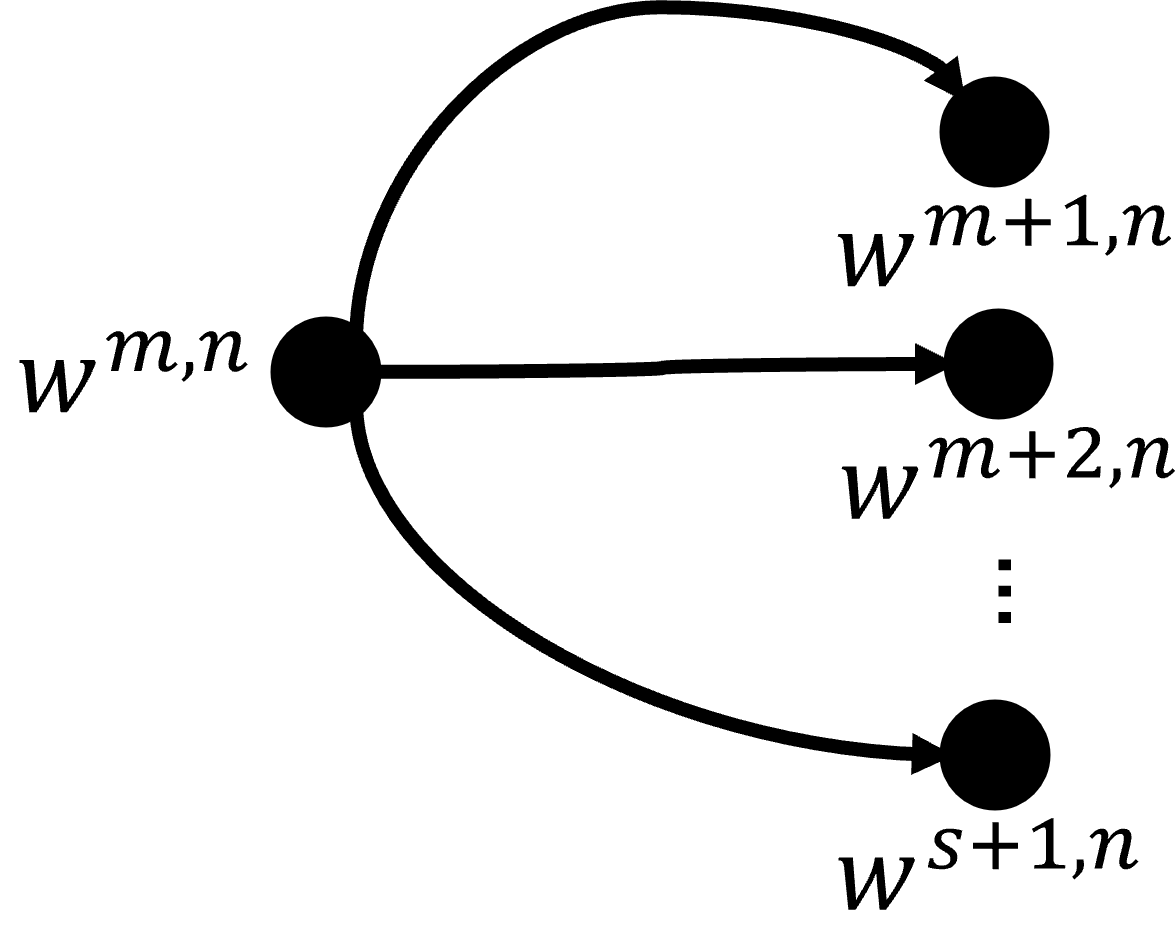}
  \caption{Corresponds to equation \ref{eq::adj2}}
  \label{fig:sfig2}
\end{subfigure}
\begin{subfigure}{.33\textwidth}
  \centering
  \includegraphics[width=.8\linewidth]{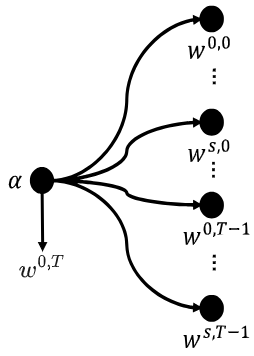}
  \caption{Corresponds to equation \ref{eq::adj3}}
  \label{fig:sfig3}
\end{subfigure}
\caption{\label{fig::adjointvariablesdefinition} Simplified computational graphs for the intermediate states and parameters. Each node corresponds to an intermediate state, and a directed edge signifies that the child intermediate state depends explicitly on the parent intermediate state. To compute the adjoint variables of the parent nodes, we require the adjoint variables of the children nodes.}
\end{figure*}

We have now established a means to compute the sensitivities of an objective function $\psi_i$, all that we are missing is a way to effectively compute the adjoint variables $\bar{w}^{0,0}_{ij}$ and $\bar{\alpha}_{ik}$. The advantage of representing the algorithm as a computational graph is that the adjoint variable of a given node $(m,n)$ can intuitively be expressed as a function of the adjoint variables of children nodes just by examining the graph. Indeed, each adjoint variable can be expressed as a summation of its children adjoint variables weighted by a partial derivative of the child node's variable with respect to the parent node's variable. To compute all the adjoint variables, we need to traverse the computational graph in reverse topological order and compute the adjoint relative to each node (with respect to a given output $i$) as a function of the adjoints of the children nodes. The recipe to compute the adjoint variables for every node are given by equations \ref{eq::adj1},\ref{eq::adj2} and \ref{eq::adj3}: 

\begin{align}
    \bar{w}^{s+1,n-1}_{ij} &= \bar{w}^{0,n}_{il}\frac{\partial w^{0,n}_l}{\partial w^{s+1,n-1}_j}\label{eq::adj1}
    \\
    \bar{w}^{m,n}_{ij}&=\sum_{r=m+1}^{s+1}\bar{w}^{r,n}_{il}\frac{\partial w^{r,n}_l}{\partial w^{m,n}_j}\label{eq::adj2}
    \\
    \bar{\alpha}_{ik} &= \frac{\partial E_i}{\partial \alpha_k}+\sum_{n=0}^{T-1}\sum_{r=0}^{s+1} 
    \bar{w}_{ij}^{r,n}\frac{\partial w^{r,n}_j}{\partial \alpha_k}\label{eq::adj3},
\end{align}
where the sub-index $i$ ranges from $1$ to $M$, the sub-indexes $j,l$ range from $1$ to $N$, the sub-index $k$ ranges from $1$ to $P$ and the computational graphs from which these equations can be induced from by inspection can be found in figure \ref{fig::adjointvariablesdefinition}. In the case that there is an explicit dependence of $E_i$ with the initial conditions, then the expression for $\bar{w}^{0,0}_{ij}$ will exceptionally be:

\begin{equation}
\bar{w}^{0,0}_{ij}=\bar{w}^{0,T}_{il}\frac{\partial E_l}{\partial u^0_j}+\sum_{r=1}^{s+1}\bar{w}^{r,0}_{il}\frac{\partial w^{r,0}_l}{\partial w^{0,0}_j}
\end{equation}

\subsubsection{Sensitivities for the trajectory cost}
Here, we consider an objective function with no end-point cost:

\begin{equation}
    \psi_i^{(2)}(\alpha) = \int_{t_0}^{t_f} R_i(\mathbf{u}(t),\alpha,t) dt \approx q_i(t_f)
\end{equation}

We use the same Runge-Kutta scheme to discretize the augmented ODE system with the $q_i$ variable. The scheme for the augmented system is: 

\begin{align}
    \mathbf{u}^{n+1} &= \mathbf{u}^n+\Delta t^n\sum_{m=1}^s b_m F(\mathbf{u}^{mm},\alpha,t^{mn}),
    \\
    q_i^{n+1} &= q_i^n+\Delta t^n\sum_{m=1}^s b_m R_i(\mathbf{u}^{mn},\alpha,t^{mn}),
\end{align}
and thus we define the following set of additional intermediate variables $v_i^{m,n}$:

\begin{align*}
    v_i^{0,0} &= q_i^0
    \\
    v_i^{m,n} &= R_i\big(\mathbf{u}^{mn},\alpha,t^{mn}\big) \qquad \text{if $m=1,\dots,s$}.
     \\
     v_i^{s+1,n} &= v_i^{0,n} + \Delta t^n \sum_{l=1}^s b_l v_i^{l,n},
     \\
     v^{0,n+1}_i &= v^{s+1,n}_i, \qquad \text{if $n>0$ and $n<T-1$},
\end{align*}
The last step is to compute the objective function, which is simply:

\begin{equation}
    v^{0,T}_i = q_i^T
\end{equation}
Similarly to before, to each intermediate variable $v_i^{m,n}$ and to each output vector entry $w^{0,T}_i$ corresponds an adjoint variable $\bar{v}_{i}^{m,n}$. We define the adjoint variables as the following total derivatives:

\begin{equation}
\bar{v}^{m,n}_{i} \defeq \frac{d \psi_i^{(2)}}{d v^{m,n}_i},
\end{equation}
exceptionally with no implied sum in the sub-index $i$. The seed will be $\bar{v}_i^{0,T}=1$ and the update equations will for the $\bar{v_i}^{m,n}$ adjoint variables will be:
\begin{align}
    \bar{v}^{s+1,n-1}_{i} &= \bar{v}^{0,n}_{i}
    \\
    \bar{v}^{m,n}_{i}&=\sum_{r=m+1}^{s+1}\bar{v}^{r,n}_{i}\frac{\partial v^{r,n}_i}{\partial v^{m,n}_i}
\end{align}
with no implied sum in the sub-index $i$. The update equations for the adjoint variables $\bar{w}^{mn}_{ij}$ now depend on $\bar{v_i}^{m,n}$: 

\begin{align}
    \bar{w}^{m,n}_{ij}&=\sum_{r=m+1}^{s+1}\bar{w}^{r,n}_{il}\frac{\partial w^{r,n}_l}{\partial w^{m,n}_j} + \bar{v}^{r,n}_{i}\frac{\partial v^{r,n}_i}{\partial w^{m,n}_j},
    \\
    \bar{\alpha}_{ik} &= \sum_{n=0}^{T-1}\sum_{r=0}^{s+1} 
    \bar{w}_{ij}^{r,n}\frac{\partial w^{r,n}_j}{\partial \alpha_k}+\bar{v}^{r,n}_{i}\frac{\partial v^{r,n}_i}{\partial \alpha_k},
\end{align}
with no implied sum in the sub-index $i$.

\subsubsection{Specifying the partial derivatives}
Given a forward pass algorithm for an explicit $s$-stage Runge-Kutta method, we can write the algorithm for its respective \textbf{reverse pass} following the adjoint update equations presented in the last section. This algorithm takes $\frac{\partial \psi_i}{\partial u^T_j}$ and $\frac{\partial \psi_i}{\partial \alpha_k}$ as inputs and returns $\bar{w}^{0,0}_{ij}$ and $\bar{\alpha}_{ik}$ as outputs. All that is left to complete the specification of this algorithm is to manually derive the partial derivatives appearing in the adjoint equations for a general $s$-stage Runge-Kutta method. The result is the following:

\begin{align}
    \frac{\partial w^{s+1,n}_i}{\partial w^{q,n}_j} &=
	\begin{cases}
		\delta_{ij}
		& \text{if } q=0\\
		\\
		\delta_{ij}\Delta t^n b_{q}
		& \text{if } 0<q<m
	\end{cases},
 \label{eq::partialderivatives1}
 \\
	\frac{\partial w^{m,n}_i}{\partial w^{q,n}_j} &=
	\begin{cases}
		J^x_{ij}\left(
		\mathbf{x}^{m,n},t^{m,n}\right)
		& \text{if } q=0\\
		\\
		J^x_{ij}\left(
		\mathbf{x}^{m,n},t^{m,n}\right)
		\Delta t^n a_{mq}
		& \text{if } 0<q<m
	\end{cases},
\label{eq::partialderivatives2}
\\
        \frac{\partial w_i^{m,n}}{\partial \alpha_k} &= \begin{cases}
		J^\alpha_{ik}(\mathbf{x}^{m,n},t^{m,n}) & \text{if } 0<m<s+1\\
            \\
		0 & \text{otherwise } 
	\end{cases},
\label{eq::partialderivatives3}
\\
	\frac{\partial w^{0,n}_i}{\partial w^{s+1,n-1}_j} &= \delta_{ij}.
\label{eq::partialderivatives4}
\end{align}
where the calculations and the definition of the Jacobian matrices $J^x$ and $J^\alpha$ and of the intermediate Runge-Kutta state $\mathbf{u^{m,n}}$ can be found in appendix \ref{appendix::calculations}. The vector-Jacobian products $\bar{w}^{m,n}_{il}J^{x,\alpha}_{lj}\left(\mathbf{u}^{m,n},\alpha,t^{m,n}\right)$ are \textbf{efficiently and automatically} evaluated using the automatic adjoint differentiation library AADC \cite{Matlogica}.  In  algorithm \ref{alg::reverse_pass}, we present the pseudo-code for the reverse pass for a single index $i$.

\begin{algorithm}
\caption{Reverse Pass}
\SetKwInOut{Input}{input}
\Input{$\frac{\partial \psi_i}{\partial u^T_j}$ and $\frac{\partial \psi_i}{\partial \alpha_k}$}
\SetKwBlock{Beginn}{beginn}{ende}
    \Begin{
    Initialise seed $\bar{w}_{ij}^{s+1,T-1} =\frac{\partial \psi_i}{\partial u_j(t_f)} $ and $\bar{\alpha}_{ik}$\;
    \For{\texttt{<$n=T-1$;$n\ge0$;$n--$>}}{
        $\bar{w}^{m,n}_{ij} =0, \text{if } 0\le m < s+1$\;
        \For{\texttt{<$m=s+1$;$m>0$;$m--$>}}{
            \For{\texttt{<$0\leq q\leq m-1$>}}{
                $\bar{w}^{q,n}_{ij}+=\bar{w}^{m,n}_{ir}\frac{\partial w^{m,n}_r}{\partial w^{q,n}_j}$
            }
        $\bar{\alpha}_{ik}+=\bar{w}^{m,n}_{ir}\frac{\partial w^{m,n}_r}{\partial \alpha_k}$
        }
    \If{$n>0$}{
    $\bar{w}^{s+1,n-1}_{ij}=\bar{w}^{0,n}_{il}\frac{\partial w^{0,n}_l}{\partial w_j^{s+1,n-1}}$}
    }
    }%\EndFor
\Return $\bar{\alpha}_{ik},\bar{w}^{0,0}_{ij}$
\label{alg::reverse_pass}
\end{algorithm}

\subsection{Checkpointing\label{sec::caching}}
The reverse pass requires knowledge of the forward trajectory whenever the Jacobian matrices depend explicitly on the state vector or the parameters. Effectively, we require the knowledge of the intermediate states $\mathbf{w}^{m,n}_i$ at the time instances $t^{m,n}$ determined adaptively at run-time. 

Obtaining the intermediate states $\mathbf{w}^{m,n}_i$ can be done in a mix of two distinct ways: either cache the intermediate states $\mathbf{w}^{m,n}_i$ and time instances $t^{m,n}$ during the forward pass or recompute the intermediate states by solving the ODE system backwards, a trade-off between performance and memory load. In their implementation of the continuous adjoint sensitivity analysis method \cite{DENSERKS2009}, the authors take the second approach by solving the ODE system backwards alongside the augmented adjoint system and use the dense output mechanism of Runge-Kutta methods to obtain the required intermediate states $\mathbf{u}^{m,n}$ at the correct instances of time, which need to be cached. In the unpublished paper \cite{zhuang2020ordinary}, the authors decided to employ a similar approach and reconstruct $\mathbf{u}^{n}$ from $\mathbf{u}^{n+1}$ by re-building the computational graph at time instances $t^{m,n}$, again requiring caching only $t^n$. This approach effectively requires the solution of the following system of equations:

\begin{equation}
    \begin{cases}
    \mathbf{w}^{s+1,n} &= \mathbf{w}^{0,n}+\Delta t^n \sum_{l=1}^s b_l \mathbf{w}^{l,n}
    \\
    \vdots
    \\
    \mathbf{w}^{m,n} &= F(\mathbf{w}^{0,n}+\Delta t^n\sum_{l=1}^{m-1}a_{mj}\mathbf{w}^{m,n},\alpha,t^{m,n})
    \end{cases}
\end{equation}
to find $\mathbf{w}^{0,n}$, which may be computationally cumbersome because of the inversion of $F$, and depending on the nature of the function might not even be possible to do analytically, thus introducing numerical error. In \cite{DENSERKS2009}, the authors need not to force the solver to compute at time instances $t^{m,n}$, whereas the time instances at which the state vectors are recomputed in \cite{zhuang2020ordinary} are determined by the forward run. However, the dense output estimation of the intermediate states in \cite{DENSERKS2009} is not exact up to machine precision. In FATODE's implementation \cite{fatode2014}, the time instances $t^n$ and every intermediate state $\mathbf{u}^{m,n}_i$ is cached, with memory cost $\mathcal{O}(N\times T\times s)$. In our implementation, we strike a balance between extra computation and memory load by caching only the hidden states $\mathbf{u}^n$ and time instances $t^n$ during the forward run, and using them to recompute the intermediate states $\mathbf{w}^{m,n}$ again, as if it were a forward run, with memory cost of $\mathcal{O}(N\times T)$.

\subsection{Vectorization via Advanced Vector Extensions}
This implementation of discrete adjoint sensitivity analysis is enhanced by utilising SIMD (Single Instruction, Multiple Data) vectorization through Advanced Vector Extensions (AVX) to explore the data-level parallelism inherent to this method. In our work, we employ AVX-256, specifically on an Intel Core i7-8850H CPU @ 2.60GHz × 12 processor, where each YMM register can hold up to four 64-bit double-precision floating point numbers.

The automatic differentiation framework of \cite{Matlogica} supports this SIMD vectorization, allowing us to efficiently compute vector-Jacobian products for multiple vectors concurrently. By leveraging AVX, we load the YMM registers, with entries from four different vectors of adjoint variables. For example, for a given set of indices $j$,$m$ and $n$, we load the YMM register with the entries $\bar{w}^{m,n}_{1j}$, $\bar{w}^{m,n}_{2j}$, $\bar{w}^{m,n}_{3j}$ and $\bar{w}^{m,n}_{4j}$, for given $j,m,n$. This functionality effectively allows us to evaluate the vector-Jacobian products of four distinct entries $i$ of the  objective function $\psi_i$ concurrently. This parallelism is particularly beneficial when there are multiple outputs of the objective function. This approach results in a performance increase when performing sensitivity analysis on a model with many parameters but few output objective functions. 

Notably, not only the vector-Jacobian products can be performed concurrently, bu the adjoint variables respective to each output can themselves be updated in parallel at a higher level of parallelism, i.e., multi-threading. Although this feature was not explored in the current work, it presents a promising direction for future research.

\section{Results and Validation \label{sec:results}}
In this section we validate our implementation of the discrete adjoint sensitivity method and test its performance compared to other sensitivity analysis software. We employ three different models for this purpose:

\begin{itemize}
    \item [1)] An $N_p\times N_p$ finite differences discretization of a two-dimensional heat equation with constant thermal diffusivity.
    \item  [2)] A Van der Pol oscillator
    \item [2)] The generalised Lotka-Volterra (GLV) equations.
\end{itemize}
We validate the algorithm by comparing the output sensitivities of several methods of sensitivity analysis with the analytical expression for the sensitivities of the heat equation. We further validate by exploring how the sensitivity error depends on the solution error for several implementations of adjoint sensitivity analysis. We provide a performance analysis by constructing work precision diagrams for our implementation and for other established sensitivity analysis software when applied to the computation of the sensitivities of the GLV equations.

The computations were performed on a laptop computer with an Intel Core i7-8850H CPU @ 2.60GHz × 12. The laptop had a 32.0 GB RAM. We ran the simulations on the Ubuntu 22.04.1 LTS operating system. We used \texttt{gcc} for compiling the C++ code and \texttt{gfortran} for compiling Fortran, both version 11.3.0, with the $-O3$ optimisation flag and we used double-precision floating point arithmetic. PETSc was installed with the following configuration:

\begin{lstlisting}[language=bash]
./config/examples/arch-ci-linux-gcc-pkgs-opt.py \
    --download-ssl --download-fblaslapack=1
\end{lstlisting}
In this section, we refer to our implementation of the discrete adjoint sensitivity analysis method simply as the "adjoint method".

\subsection{Two-dimensional heat equation with classic RK methods}
Consider the two-dimensional heat equation in the domain $\Omega = [0,1]\times[0,1]$ and time interval $t\in(0,t_f)$ with homogeneous Dirichlet boundary conditions:

\begin{equation}
   \frac{\partial u(\mathbf{x},t)}{\partial t} = \alpha \nabla u(\mathbf{x},t), \qquad
    u(\mathbf{x},0) = u^0(\mathbf{x}), \qquad
    u(\mathbf{x},t) = 0, \forall \mathbf{x}\in\partial\Omega,
    \label{eq::heat_equation_2D}
\end{equation}
where $u(\mathbf{x},t)$ is the temperature distribution in the domain $\Omega$ and $\alpha$ is the constant thermal diffusivity. The solution is given by,

\begin{equation}
    u(\mathbf{x},t) = \sum_{m=1}^\infty\sum_{n=1}^\infty \gamma_{mn}\sin(m\pi x)\sin(n\pi y)\exp{(-\lambda_{mn}^2 t)}
\end{equation}
where $\lambda_{mn}$ is the normal mode eigenvalue and $\gamma_{mn}$ the 2D Fourier coefficient of the initial condition $u^0(\mathbf{x})$,

\begin{equation}
    \lambda_{mn}= \sqrt{\alpha\pi^2(m^2+n^2)},
    \qquad
    \gamma_{mn} = 4\int_{0}^{1}\int_{0}^{1} u(\mathbf{x},0)\sin(m\pi x)\sin(n\pi x) \text{ }dxdy.
\end{equation}
The sensitivity of the solution $u(\mathbf{x},t_f)$ at the final time with respect to the constant thermal diffusivity, $\Bar{\alpha}(\mathbf{x})\defeq\frac{\partial u(\mathbf{x},t_f)}{\partial \alpha}$, is given in terms of the closed-form analytic expression:

\begin{equation}
    \Bar{\alpha}(\mathbf{x}) = -\pi^2 t_f\sum_{m,n=1}^\infty (m^2+n^2)\gamma_{mn}\times\sin(m\pi x)\sin(n\pi y)\exp{(-\lambda_{mn}^2 t_f)}.
\label{exactadjointheat}  
\end{equation}

We discretize equation \ref{eq::heat_equation_2D} in space using a finite-difference scheme. We define an evenly spaced grid $(x_i,y_j)$ of $N_p\times N_p$ points, where $x_i = (i-1)\Delta x$ and $y_j = (j-1)\Delta y$, with $i,j = 1,\dots,N_p$. Let $u_{k}\approx u(x_i,y_j)$, with $k = i+N_p(j-1)$, and define the set of interior grid points as $I = \left\{k\in\mathbb{N}\colon (x_{i(k)},y_{j(k)}) \in \text{int}(\Omega)\right\}$ and the set of boundary points as $B = \left\{k\in\mathbb{N}\colon (x_{i(k)},y_{j(k)}) \in \text{fr}(\Omega)\right\}$. Equation \ref{eq::heat_equation_2D} is thus approximated by the following system of ODEs:

\begin{align}
    \frac{d u_{k}(t)}{dt} &=
    \alpha \bigg(
    \frac{u_{k-1}-2u_{k}+u_{k+1}}{\Delta x^2}\nonumber
    +\frac{u_{k-N_p}-2u_{k}+u_{k+N_p}}{\Delta y^2}
    \bigg), \forall k \in I,\nonumber
    \\
    \frac{du_{k}(t)}{dt} &= 0\qquad \forall k \in B, \qquad u_k(0) = u^0(x_{i(k)},y_{j(k)}) \qquad \forall k \in I\cup B.
\label{eq::heat_equation_ode_system}
\end{align}
In the following experiments, we numerically compute  the sensitivities with respect to the thermal diffusivity with several methods of sensitivity analysis and compare them with the analytical solution. We consider $\Delta y = \Delta x = \frac{1}{N_p-1}$, $u^0(\mathbf{x}) = \sin(\pi x)\sin(\pi y)$ and $\alpha = 1$. The ODE system is integrated in the interval $t \in [0,10^{-2}]$ with a constant time step $\Delta t$ using either the explicit Euler solver or the classic $4$th order Runge-Kutta solver provided in the \texttt{boost::odeint} library.

\subsubsection{Results and Discussion}
In tables \ref{tab:table1}, \ref{tab:table2} and \ref{tab:table3} we present the relative error of the forward solution along with the relative error for the sensitivities computed using different sensitivity analysis methods across various grid sizes $N_p$ \footnote{We remark that the implementation of all these methods has been done from scratch by the authors of this article.}.
The forward solvers used to obtain the results were: the explicit Euler method with $\Delta t = 5\times10^{-5}$ for table \ref{tab:table1} and the classic $4$th order Runge-Kutta (RK4) method with $\Delta t = 5\times10^{-5}$ and $\Delta t = 1\times10^{-5}$ for tables \ref{tab:table2} and \ref{tab:table3} respectively. The time steps were chosen based on the stability condition $\Delta t<\frac{1}{4\alpha} \Delta x^2$ for the largest grid size $N_p$

We define the relative error $\epsilon_r$ as:

\begin{equation}
    \epsilon_r = \frac{||\Bar{\alpha}_i-\Bar{\alpha}(\mathbf{x}_i)||_{\infty}}{||\Bar{\alpha}(\mathbf{x}_i)||_{\infty}}
\end{equation}
where $\Bar{\alpha}_i$ is the numerical solution computed via a sensitivity analysis algorithm and $\Bar{\alpha}(\mathbf{x}_i)$ is the exact solution evaluated at the grid points $\mathbf{x}_i$. The error for the forward solution is defined analogously. 

\begin{table}[h!]
\caption{\label{tab:table1} Relative errors (\%) for the Euler method with $dt = 5\times 10 ^{-5}$ and bump $\Delta \alpha=10^{-9}$
}
\centering
\begin{tabular}{|p{0.4cm}|p{1.5cm}|p{1.5cm}|p{1.5cm}|p{1.4cm}|p{1.4cm}|}
    \hline
    \multirow{2}{*}{$N_p$} & \multicolumn{5}{|c|}{Relative Error ($\times 10^{-1}$ \%)} \\
    \cline{2-6}
                             &\multicolumn{1}{|c|}{Forward}
                             &\multicolumn{1}{|c|}{Adjoint} 
                             & \multicolumn{1}{|c|}{CFSA}
                             & \multicolumn{1}{|c|}{ND}
                             & \multicolumn{1}{|c|}{CASA}
                             \\
                             \cline{2-6}
    \hline
    10   & $8.827$&$7.260$&   $7.260$ & $7.260$& $3.366$\\
    \hline
    30    &$1.392$ &$0.104$&   $0.104$ & $0.105$& unstable\\
    \hline
    50    &$0.731$ &$0.615$&   $0.615$ & $0.615$& unstable\\
    \hline
\end{tabular}
\end{table}

\begin{table}[h]
\caption{\label{tab:table2} Relative errors (\%) for classic RK4 method with $dt = 5\times 10 ^{-5}$ and bump $\Delta \alpha=10^{-9}$
}
\centering
\begin{tabular}{|p{0.4cm}|p{1.5cm}|p{1.5cm}|p{1.5cm}|p{1.4cm}|p{1.4cm}|}
    \hline
    \multirow{2}{*}{$N_p$} & \multicolumn{5}{|c|}{Relative Error ($\times 10^{-1}$\%)} \\
    \cline{2-6}
                             &\multicolumn{1}{|c|}{Forward}
                             &\multicolumn{1}{|c|}{Adjoint} 
                             & \multicolumn{1}{|c|}{CFSA}
                             & \multicolumn{1}{|c|}{ND}
                             & \multicolumn{1}{|c|}{CASA}
                             \\
                             \cline{2-6}
    \hline
    10   & $9.271$&$7.981$&   $8.135$ & $8.137$& $3.175$\\
    \hline
    30   & $2.807$ &$0.629$&   $0.785$ & $0.789$& unstable\\
    \hline
    50   & $1.658$    &$0.119$&   $0.275$ & $0.275$& unstable\\
    \hline
\end{tabular}
\end{table}

\begin{table}[h!]
\caption{\label{tab:table3} Relative errors (\%) for classic RK4 method with $dt = 1\times 10 ^{-5}$ and bump $\Delta \alpha=10^{-9}$
}
\centering
\begin{tabular}{|p{0.4cm}|p{1.5cm}|p{1.5cm}|p{1.5cm}|p{1.4cm}|p{1.4cm}|}
    \hline
    \multirow{2}{*}{$N_p$} & \multicolumn{5}{|c|}{Relative Error ($\times 10^{-1}$\%)} \\
    \cline{2-6}
                             &\multicolumn{1}{|c|}{Forward}
                             &\multicolumn{1}{|c|}{Adjoint} 
                             & \multicolumn{1}{|c|}{CFSA}
                             & \multicolumn{1}{|c|}{ND}
                             & \multicolumn{1}{|c|}{CASA}
                             \\
                             \cline{2-6}
    \hline
    10   & $9.271$& $8.104$&   $8.135$ & $8.155$& $7.639$\\
    \hline
    30   & $2.807$ &$0.754$&   $0.785$ & $0.785$& unstable\\
    \hline
    50   & $2.807$ & $0.245$&   $0.275$ & $0.274$& unstable\\
    \hline
    \end{tabular}
\end{table}

The discrepancies between the analytical and numerical sensitivities can likely be attributed to the discretization error associated with the approximation of the heat equation by the system of ODEs, as well as the error introduced by the explicit numerical method used for discretization due to a finite time step $\Delta t$. Additionally, the numerical differentiation method is also affected by truncation error and cancellation error. 

Despite these factors, all the methods perform similarly in terms of accuracy for the simulation parameters considered. Our implementation of the discrete adjoint sensitivity analysis method yields a negligibly better result in most cases, whereas the CASA method appears more accurate when stable. This instability is likely due to inaccuracies in reconstructing the forward trajectory during the reverse pass in this particular implementation of the method. Notably, the results when using the Euler method are better than those of the RK4 method for $N=10$ and $N=30$, but not for $N=50$, and we notice that increasing the spatial resolution too much degrades the accuracy of the Euler method, a phenomenon not observed with the RK4 method. Additionally, decreasing the time step does not improve the results for the RK4 method. 

In conclusion, our implementation of the discrete adjoint sensitivity analysis method performs correctly and provides results comparable to other sensitivity analysis methods.

\subsection{Van der Pol oscillator with adaptive RK methods}
In this section we consider the Van der Pol oscillator, a second-order non-linear ODE with the following dynamics:

\begin{equation}
        \frac{d^2x}{dt^2} = \mu(1-x^2)\frac{dx}{dt}-\mu x,
\end{equation}
which can be written as the following system of first-order ODEs,
\begin{equation}
    \frac{d}{dt}
    \begin{bmatrix}
    x\\
    v
    \end{bmatrix}
    =
    \begin{bmatrix}
    v\\
    \mu(1-x^2)v-\mu x
    \end{bmatrix}.
\end{equation}
PETSc TSAdjoint provides an example of sensitivity analysis for this system in the tutorial file "\textit{tutorials/ex20adj.c}".

\begin{figure}[h!]
    \centering
    \includegraphics[width = 0.65\textwidth]{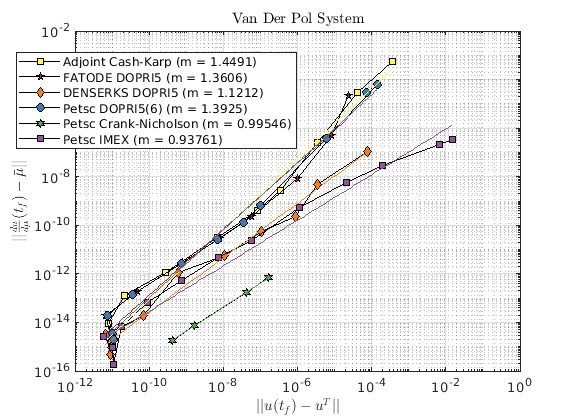}
    \caption{Plot of the solution error ($x$-axis) against the sensitivity error ($y$-axis) for the Van der Pol equations for multiple implementations of adjoint sensitivity analysis.}
    \label{fig:vanderpol}
\end{figure}

We solved the system over the interval $t \in [0,5\times10^{-1}]$ with a parameter value $\mu = 1\times 10^3$ and initial conditions $x(0) = 2$ and $v(0) = -\frac{2}{3} + \frac{10}{81\mu}- \frac{292}{2187\mu^2}$. The system was solved using the same absolute and relative error tolerance for the forward solver, ranging from $10^{-3}$ to $10^{-15}$. We computed the sensitivities $\frac{dx}{d\mu}$ and $\frac{dv}{d\mu}$ using several state-of-the-art implementations of adjoint sensitivity analysis and compared the results to our implementation. We considered different forward solvers depending on the implementation of adjoint sensitivity analysis. As for explicit adaptive methods, the Cash-Karp method was used in our implementation while the Dormand-Prince method was used for FATODE, DENSERKS and PETSc TSAdjoint. Both these ERK methods are of the same order and number of stages. Additionally, we also considered the Crank-Nicholson method and an Implicit-Explicit (IMEX) method in PETSc TSAdjoint to provide an additional baseline. We plotted the solution error against the sensitivity error, computed as the inf-norm of the difference between the solution of each implementation and a finely resolved solution with an absolute and relative error tolerance of $2.5\times 10^{-14}$, obtained using MATLAB's recent feature for sensitivity analysis. Results are presented in Figure \ref{fig:vanderpol}.

\subsubsection{Results and Discussion}
Notably, the PETSc methods exhibit better sensitivity results overall for this problem, with this advantage becoming more pronounced as error tolerances increase. The sensitivity error is approximately one order of magnitude smaller to our method, FATODE, and DENSERKS when the solution error is small, but the gap increases for higher solution error. Our method, FATODE, and DENSERKS exhibit similar results across different solution errors. However, our method shows the highest sensitivity errors, happening for the highest error tolerances. This result could be attributed to our use of the Cash-Karp method instead of the Dormand-Prince method. All graphs show a cut-off around $10^{-11}$ , due to the solutions being compared against a finite-precision MATLAB baseline solution.

In summary, all methods perform comparably. However, the PETSc methods show slightly better accuracy for this specific problem. This example provides additional confidence in the correctness of our method.

\subsection{Generalized Lotka-Volterra Equations with adaptive RK methods}
The generalized Lotka-Volterra (GLV) equations form a system of ODEs that model the dynamics of the populations of $N$ species $\{x_i\}_{i=1}^N$:

\begin{equation}
    \frac{d x_i(t)}{dt} = x_i(r_i+A_{ij}x_j(t)), \qquad x_i(0) = x_i^0 \qquad (\text{no sum in $i$})
\end{equation}
where $\mathbf{r}\in\mathbb{R}^N$ is a vector of intrinsic growth/death rates, representing the rate of change of the populations in the absence of interactions, and $A\in\mathbb{M}_{N\times N}(\mathbb{R})$ is a matrix of interaction coefficients, quantifying the effect that each population has on the dynamics of others. The original Lotka-Volterra equations were first studied independently by Alfred J. Lotka in \cite{lotka2002contribution} (reprint) and Vito Volterra in \cite{volterra1926variazioni}.

Consider the following parameters vector $\alpha$:

\begin{equation}
    \alpha_{k(m,n)} = 
    \begin{cases}
        r_n, \qquad &\text{if } m = 0\\
        A_{mn},  \qquad &\text{if } m > 0
    \end{cases},
\end{equation}
where $m= 0,\dots,N$ and $n=1,\dots,N$ and by defining the mapping $k(m,n)$ between parameter vector indexes $k$ and the indexes of the vector of growth rates and interaction matrix as $k(m,n)=m\times N+n$. In doing so, it is straightforward to cast the GLV system in the form of system \ref{eq::system_of_odes}:

\begin{align}
    \frac{\partial x_i}{\partial t} &= F_i(\mathbf{x},\mathbf{\alpha}), \qquad x_i(0) = x_i^0 \nonumber
    \qquad 
    \text{with}
    \qquad 
    \\
    \qquad 
    F_i(\mathbf{x},\alpha) &= x_i(\alpha_{k(0,i)}+\sum_{j=1}^N\alpha_{k(i,j)}x_j) \qquad \text{(no sum in $i$)}
    \label{eq::GLV_system}
\end{align}
The above is a system of $N$ equations with up to $N+N^2$ parameters, encompassing both the growth rates and interaction matrix coefficients. This complexity makes the GLV system an ideal candidate for evaluating adjoint sensitivity methods, particularly due to the significantly larger number of parameters compared to the number of objective functions.

In this study, GLV systems of varying sizes $N$ are simulated using randomly generated interaction matrices $A_{ij}$ to ensure stable solutions, following Robert M. May's approach \cite{may1972will}. We considered the growth rates $r_i = 0.1$, initial conditions $x_i^0 = 0.1$ and simulate over the time interval $t\in(0,10)$. The sensitivity matrix $\frac{d u_i^T}{d\alpha_k}$ is computed using several methods and implementations of sensitivity analysis, for which the forward solution is computed using adaptive RK methods, considering $\alpha$ as the vector encompassing all the $N+N^2$ possible parameters. The performance and efficacy of these methods and implementations are analysed and compared.

\subsubsection{Results and Discussion}
We present the results in three subsections: the first provides a performance comparison between our implementation of the discrete adjoint method and our implementation of other methods of sensitivity analysis; the second subsection provides a performance comparison between our implementation of the discrete adjoint method, FATODE's implementation of the discrete adjoint method \cite{fatode2014}, DENSERKS' implementation of the continuous adjoint method \cite{DENSERKS2009} and finally PETSc TSAdjoint's \cite{petsc2022} implementation of the discrete adjoint method; the third subsection provides a comparison between all algorithms regarding how the execution time scales with the parameter count $P$.

To generate work-precision diagrams, we compute the forward solution using a range of local absolute and relative error tolerances from $10^{-4}$ and $10^{-14}$ (the absolute and relative error tolerances are taken to be the same), and measure the global sensitivity error $\epsilon$ as the $||\cdot||_\infty$-norm of the difference between the algorithm outputs (the sensitivities with respect to the parameters for all the outputs of the ODE at $t=t_f$ for a given tolerance) and a baseline solution computed with the same approach but with a tolerance of $10^{-15}$. Our experiments show that this is qualitatively equivalent to using the CFSA method with an RK Cash-Karp5(4) solver at a tolerance of $10^{-15}$ to generate the baseline output sensitivities, though this approach becomes unfeasible for very large systems. We measure the average execution time of computing the sensitivities over multiple runs. We consider a bump of $\Delta\alpha = 1\times 10^{-6}$ for the numerical differentiation ( finite differences) method.

\begin{figure*}[h!]
        \centering
        \begin{subfigure}[]{0.475\textwidth}
            \centering
            \includegraphics[width=\textwidth]{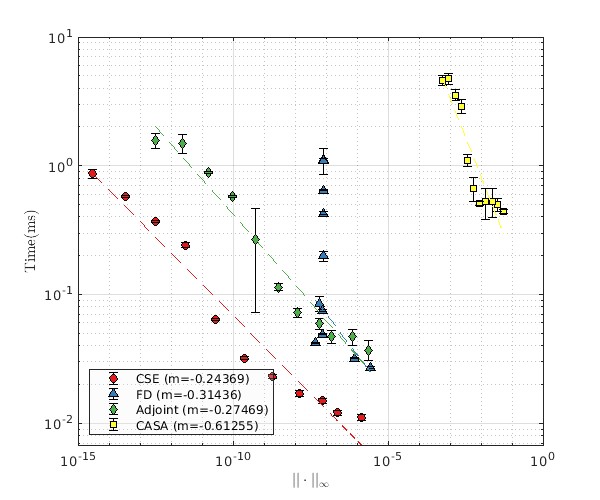}
            \caption[]%
            {{\small $N=2$}}    
            \label{fig::WorkPrecisionN2}
        \end{subfigure}
        \hfill
        \begin{subfigure}[]{0.475\textwidth}  
            \centering 
            \includegraphics[width=\textwidth]{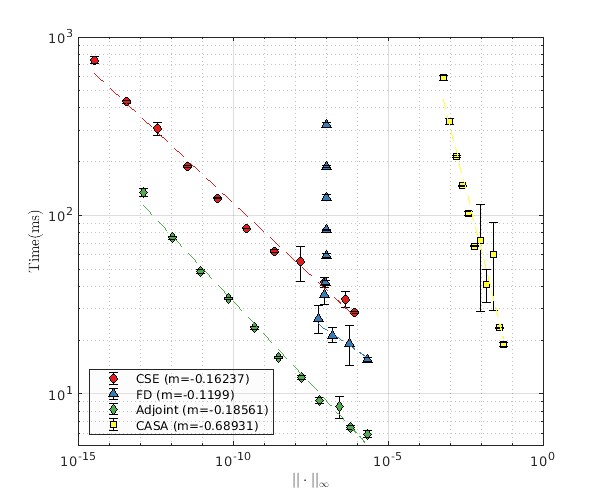}
            \caption[]%
            {{\small $N=10$}}    
            \label{fig::WorkPrecisionN10}
        \end{subfigure}
        \vskip\baselineskip
        \begin{subfigure}[]{0.475\textwidth}   
            \centering 
            \includegraphics[width=\textwidth]{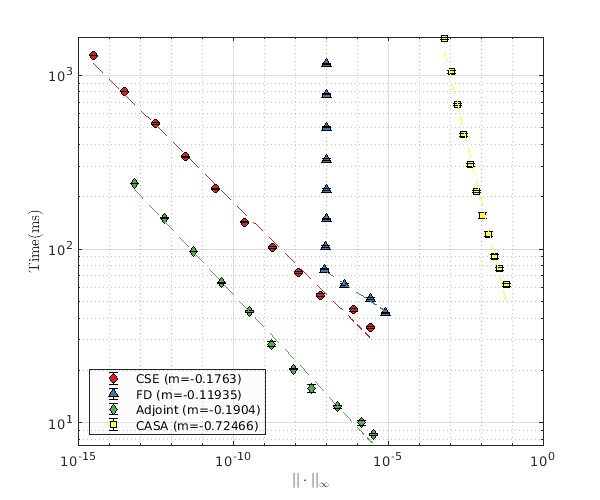}
            \caption[]%
            {{\small $N=30$}}    
            \label{fig::WorkPrecisionN30}
        \end{subfigure}
        \hfill
        \begin{subfigure}[]{0.475\textwidth}   
            \centering 
            \includegraphics[width=\textwidth]{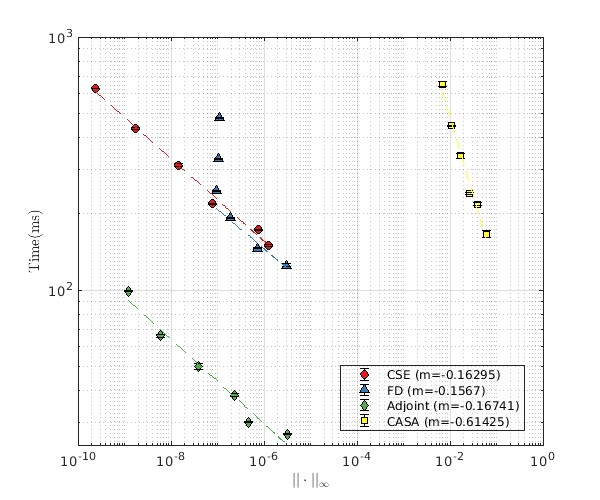}
            \caption[]%
            {{\small $N=40$}}    
            \label{fig::WorkPrecisionN40}
        \end{subfigure}
        \caption[ The average and standard deviation of critical parameters ]
        {\small Average execution clock time in milliseconds plotted against the global sensitivity error for several methods of sensitivity analysis. The average was taken over $30$ executions. The vertical lines represent a standard deviation.
        The dashed lines correspond to a linear regression of the data to the power law $t \sim \epsilon^{m}$ ignoring the data corresponding to the largest three error values.} 
\end{figure*}

\paragraph{Comparison between different methods of sensitivity analysis}
In figures \ref{fig::WorkPrecisionN2},\ref{fig::WorkPrecisionN10}, \ref{fig::WorkPrecisionN30} and \ref{fig::WorkPrecisionN40} we present work-precision diagrams for GLV systems of sizes $N=2$, $N=10$, $N=30$ and $N=40$ respectivelly. In these graphs the following methods figure:

\begin{itemize}
    \item The implementation of the discrete adjoint method figured in this article is shown in green.
    \item A handwritten implementation of the continuous forward sensitivity method (the Jacobian-vector products are manually implemented) is shown in red.
    \item A numerical differentiation method with a bump of $1\times 10^{-6}$ is shown in blue.
    \item An implementation of the continuous adjoint sensitivity method using automatic differentiation (via AADC) for the vector-Jacobian products, which does not take into account the forward trajectory when performing the reverse pass, is shown in yellow.
\end{itemize}
All the above methods use the \texttt{boost::numeric::odeint} Runge-Kutta Cash-Karp5(4) method as forward integrator.

The adjoint method performs better than the other algorithms when the number of parameters $P$ is large, as expected. Its performance becomes competitive with the CFSA method at around $N=10$, for which approximately the same global error is achieved at the same computational cost. The method becomes consistently faster compared to the others when $N$ increases past $10$, reaching almost one order of magnitude speed up compared to the CFSA method when $N=40$. 

The global error for the ND method stagnates at around $10^{-7}$ regardless of the system size $N$. This is due to the fact that the increase in accuracy resulting from an increase in the local error tolerance starts being negligible relative to catastrophic cancellation error. This phenomenon happens regardless of the bump that is considered, but the stagnation may occur at different error values depending on the bump. However, decreasing the bump further past $1\times 10^{-6}$ does not mean that the results will be better, existing an optimal bump for which the results are the best. The method is far from ideal for computing sensitivities with strict error tolerances. 

We notice that the CFSA method achieves the lowest global error of around $2\times 10^{-15}$ in all the graphs, but our hypothesis is that this happens because the error is compared against a fine CFSA solution and not a real analytical solution. This behaviour is to be expected and is in no way an advantage of this method. The method seems to deliver good results consistently, but is expected to be slower than the discrete adjoint method when $P\gg N$.

The CASA method gives the worst global errors of all methods, never falling below $10^{-4}$, even when the local error tolerance is taken to be $10^{-14}$. This is due to the fact that, in this implementation, the forward trajectory is not taken into account when performing the reverse pass.

\begin{figure*}[h!]
        \centering
        \begin{subfigure}[]{0.475\textwidth}
            \centering
            \includegraphics[width=\textwidth]{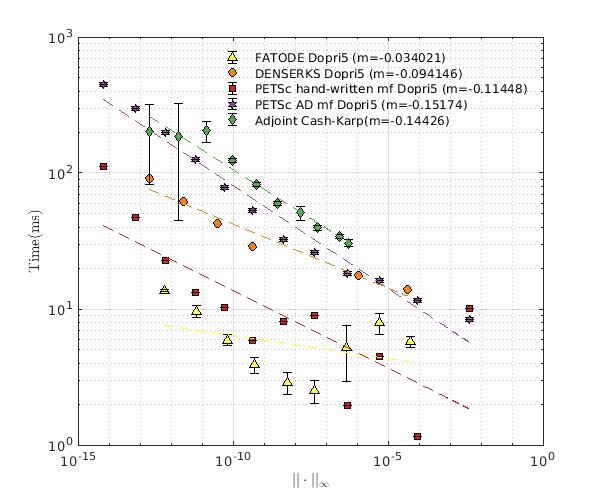}
            \caption[]%
            {{\small $N=10$}}  
            \label{fig:N10}
        \end{subfigure}
        \hfill
        \begin{subfigure}[]{0.475\textwidth}  
            \centering 
            \includegraphics[width=\textwidth]{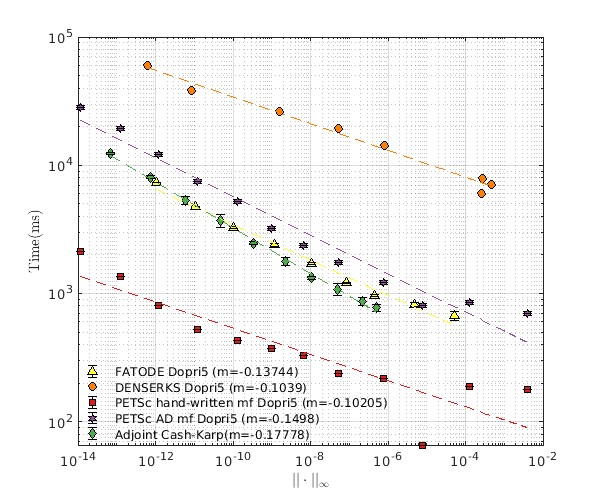}
            \caption[]%
            {{\small $N=55$}}    
            \label{fig:N55}
        \end{subfigure}
        \vskip\baselineskip
        \begin{subfigure}[]{0.475\textwidth}   
            \centering 
            \includegraphics[width=\textwidth]{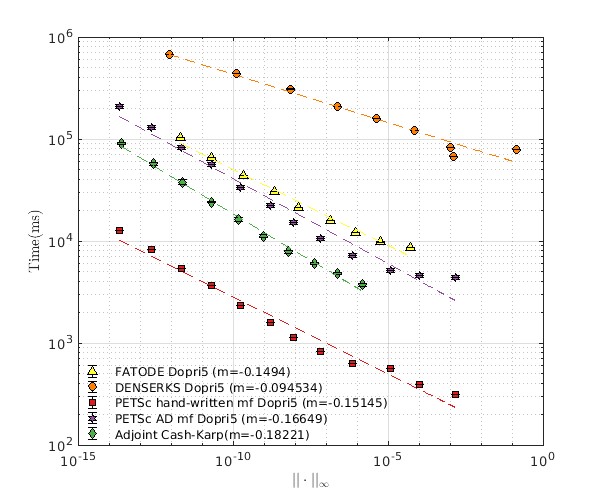}
            \caption[]%
            {{\small $N=100$}}    
            \label{fig:N100}
        \end{subfigure}
        \hfill
        \begin{subfigure}[]{0.475\textwidth}   
            \centering 
            \includegraphics[width=\textwidth]{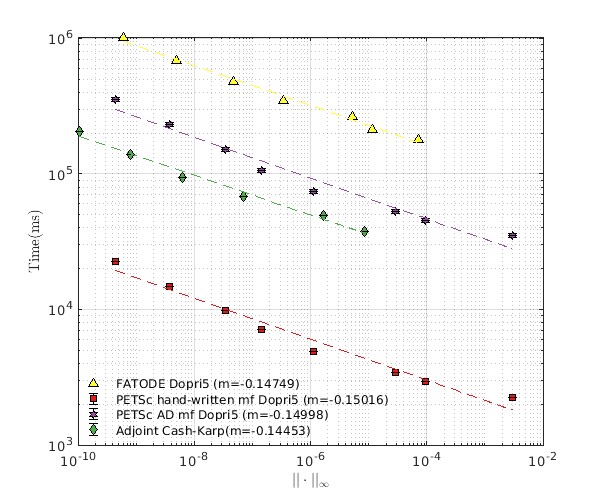}
            \caption[]%
            {{\small $N=200$}}    
            \label{fig:N200}
        \end{subfigure}
        \caption[ ]
        {\small The average execution time (in milliseconds) is plotted against the global sensitivity error for five different implementations of adjoint sensitivity analysis applied to GLV systems of varying sizes. Vertical error bars indicate one standard deviation. Dashed lines represent linear regressions of the data following the power law $t \sim \epsilon^{m}$, with the slope of each line provided in the plot legend. Data for DENSERKS at $N=200$ is absent due to excessive execution time.}
\end{figure*}

\paragraph{Comparison between different implementations of adjoint methods}
Figures \ref{fig:N10},\ref{fig:N55},\ref{fig:N100} and \ref{fig:N200} show work-precision diagrams for the sensitivity analysis of GLV systems of size $N=10$, $N=55$, $N=100$ and $N=200$. The primary distinctions between the algorithms compared are in the selection of the forward solver—either Dormand-Prince or Cash-Karp methods—and in the techniques used for evaluating vector-Jacobian products. These techniques include matrix-free approaches (using automatic differentiation or a hand-coded vector-Jacobian product) and explicit approaches (where the Jacobian is first computed and then used to evaluate the product). The following five implementations are included in the comparison:

\begin{itemize}
    \item The implementation of the discrete adjoint method figured in this article, using the Runge-Kutta Cash-Karp5 as forward solver method and automatic differentiation to evaluate the vector-Jacobian products, in green.
    \item FATODE's implementation of the discrete adjoint method using the Runge-Kutta Dopri5 as forward solver method with explicit Jacobian evaluation, in yellow.
    \item DENSERKS' implementation of the continuous adjoint method using the Runge-Kutta Dopri5 method as forward solver  and a hand-written implementation of the vector-Jacobian products, in orange.
    \item PETSc TSAdjoint's implementation of the discrete adjoint method using the Runge-Kutta Dopri5 method as forward solver in red with a hand-written implementation of the vector-Jacobian products, in red.
    \item PETSc TSAdjoint's implementation of the discrete adjoint method using the Runge-Kutta Dopri5 method as forward solver and the same automatic differentiation mechanism as our implementation to evaluate the vector-Jacobian products, in purple.
\end{itemize}

The performance comparison of the aforementioned implementations of adjoint sensitivity analysis reveals distinct behaviors based on system size and complexity. For the smallest system size, $N=10$, our algorithm is slower than all other implementations, with performance more than an order of magnitude behind FATODE and PETSc with hand-written vector-Jacobian products. Our implementation and PETSc with automatic differentiation perform similarly, though both are slightly slower than DENSERKS. At $N=55$, our implementation matches FATODE's performance, with both being slightly faster than PETSc with automatic differentiation. DENSERKS is nearly an order of magnitude slower, while PETSc with hand-written vector-Jacobian products is about an order of magnitude faster. For $N=100$ and $N=200$, the benefits of vectorization become apparent, as our implementation outperforms all others except for PETSc with hand-written vector-Jacobian products, which remains roughly an order of magnitude faster. While FATODE remains competitive at $N=100$, it falls behind at $N=200$, where our implementation and PETSc with automatic differentiation are significantly faster. Data for DENSERKS at $N=200$ is unavailable due to prohibitively long execution times.

These performance differences are likely due to specific implementation details, which we proceed to describe. 

\begin{itemize}
    \item FATODE's implementation caches Runge-Kutta intermediate states $u^{m,n}$, the state vector $u^n$, and time steps $dt^n$ in memory, avoiding re-computation and resulting in faster performance but at the cost of high memory usage. In contrast, our approach follows the strategy described in Section \ref{sec::caching}, which involves caching only the state vector and not the Runge-Kutta intermediate states. An inefficiency in FATODE’s implementation is that the Jacobian is explicitly computed before the vector-Jacobian product, whereas in our implementation and PETSc the vector-Jacobian products are directly computed via automatic differentiation, contributing to improved performance. Additionally, our implementation leverages vectorization to compute the vector-Jacobian products relative to different outputs concurrently, whereas FATODE does not employ any parallelization strategy.
    \item DENSERKS uses a two-level checkpointing scheme, storing time steps, time instances, and state vectors every $N_d$ steps, which requires additional computation to recover Runge-Kutta intermediate states. In this implementation, the vector-Jacobian product is efficiently evaluated through a hand-written routine. Notably, DENSERKS sometimes diverges at certain local error tolerances (not shown in the figures), an issue not observed in FATODE, PETSc, or our implementation of the discrete adjoint method. This may occur because DENSERKS solves the adjoint system of ODEs backwards, which is generally stiffer than the original ODE system \cite{juliaarticle}. Since the same tolerances are used for both the forward and adjoint solvers, and the adjoint system is stiffer, the solver may diverge, despite taking advantage of the dense output mechanism of explicit Runge-Kutta methods. Although DENSERKS allows for stricter tolerances for the adjoint solver, they can only be set heuristically, and setting a stricter tolerance can lead to more time steps and slower performance during the adjoint solve. In contrast, the number of time steps for the adjoint solve in FATODE, PETSc, and our implementation is the same as for the forward pass, demonstrating the advantage of the discrete adjoint method over the continuous adjoint method. DENSERKS does not employ any parallelization strategy.
    \item PETSc is able to employ an optimal checkpointing strategy that minimizes re-computations while managing memory constraints by selectively storing intermediate states, balancing memory usage and computational speed. However, in both implementations discussed here, PETSc stores the state vector $u^n$ at each time step and recomputes the Runge-Kutta intermediate states, the same approach as our implementation. For large $N$, our implementation consistently outperforms PETSc with automatic differentiation, suggesting that vectorization provides a performance advantage when computing vector-Jacobian products concurrently. Nonetheless, PETSc with hand-written vector-Jacobian products outperforms all other implementations, likely because it avoids the overhead introduced by AD, although this approach requires the user to derive an analytic expression for the vector-Jacobian product. In this implementation, no parallelization features of PETSc are utilised. 

\end{itemize}
FATODE and DENSERKS are compatible with AD through tools like TAMC \cite{Giering1998a} or TAPENADE \cite{hascoet2013tapenade}, requiring manual integration of generated code into sensitivity analysis programs. In practice, the user has to copy-paste the generated code from these tools into their own sensitivity analysis program. Our method uses AD in the operator-overloading approach via the AADC library, and the code for the derivative computation is generated during run-time and is integrated in our sensitivity analysis library. In \cite{petscautodiff} the authors describe how to integrate automatic differentiation in a PETSc implementation of sensitivity analysis using ADOL-C \cite{adol1996}. Our approach employs AADC \cite{Matlogica} to implement a matrix-free approach of the reverse pass.

%\begin{table}[ht]
%\centering
%\small
%\setlength{\tabcolsep}{6pt} % Reduce the space between columns
%\begin{tabular}{|p{2cm}|p{3cm}|p{3cm}|p{3cm}|} % Adjust the width as needed
%\hline
%\textbf{Implementation} & \textbf{Check-pointing} & \textbf{Vector-Jacobian Product} & %\textbf{Parallelisation} \\
%\hline
%Adjoint & \parbox{3cm}{\raggedright Cache only the solution in memory} & \parbox{3cm}{\raggedright Efficient} & \parbox{3cm}{\raggedright Vectorization of vector-Jacobian products via AVX} \\
%\hline
%FATODE & \parbox{3cm}{\raggedright Cache the solution, time and intermediate states in memory} & \parbox{3cm}{\raggedright Explicit Jacobian evaluation} & \parbox{3cm}{\raggedright None }\\
%\hline
%DENSERKS & \parbox{3cm}{\raggedright Two-level check-pointing} & \parbox{3cm}{\raggedright Efficient %if implemented by hand} & \parbox{3cm}{\raggedright None} \\
%\hline
%TSAdjoint & \parbox{3cm}{\raggedright Optimal check-pointing} & \parbox{3cm}{\raggedright Efficient} %& \parbox{3cm}{\raggedright Inherits parallel infrastructure of PETSc solvers.} \\
%\hline
%\end{tabular}
%\caption{Comparison of implementation details}
%\label{tab:implementation-details}
%\end{table}

\begin{figure*}[ht]
    \centering
    \begin{subfigure}[]{0.5\textwidth}  
        %\centering
        \includegraphics[width=\linewidth]{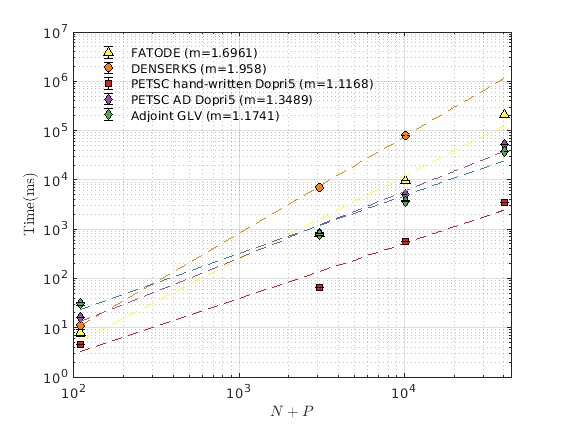}
        \caption{Local error tolerance of $10^{-5}$}
        \label{fig:fig1}
    \end{subfigure}%
    \begin{subfigure}[]{0.5\textwidth}  
        %\centering
        \includegraphics[width=\linewidth]{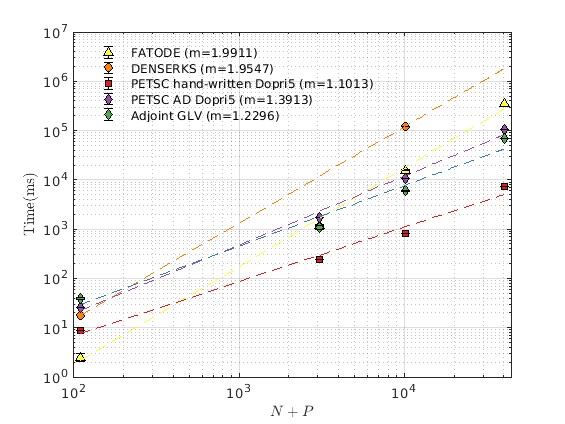}
        \caption{Local error tolerance of $10^{-7}$}
        \label{fig:fig2}
    \end{subfigure}%
    \vskip\baselineskip
    \begin{subfigure}[]{0.5\textwidth}  
        \centering
        \includegraphics[width=\linewidth]{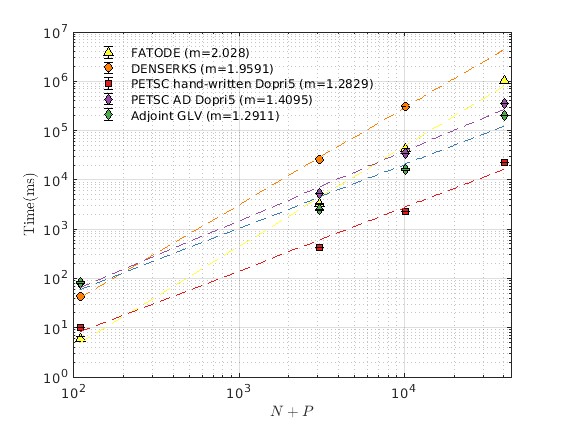}
        \caption{Local error tolerance of $10^{-10}$}
        \label{fig:fig3}
    \end{subfigure}
    \caption{\label{fig:fig}Plot of execution clock time in milliseconds against the number of parameters in log-log scale. The dashed lines correspond to linear regressions to the power law  $t \sim  P^{m}$.}
\end{figure*}

\paragraph{Dependence with number of parameters\label{paragrah::timeVSnumparameters}}
In figures \ref{fig:fig1} and \ref{fig:fig2} we present the average execution clock time in milliseconds as a function of the number of $N+P$ for given forward solver local error tolerances.

Upon reviewing figures \ref{fig:fig1}, \ref{fig:fig2} and \ref{fig:fig3}, it is evident that the computational time for the adjoint method implementation presented in this article scales similarly to that of PETSc’s hand-written implementation, regardless of local error tolerance. In these cases, the computational time scales approximately linearly with $N+P$, which aligns with expectations for adjoint sensitivity analysis methods \cite{juliaarticle}. PETSc with automatic differentiation also exhibits near-linear scaling with $N+P$, but consistently reveals slightly worst scaling of the computational time with the number of parameters than FATODE and our implementation across various local error tolerances. Finally, linear scaling does not hold for the FATODE and DENSERKS implementations. At lower error tolerances, the plots suggest that their scaling approaches quadratic in $N+P$. FATODE's deviation from linear scaling could be due to its explicit computation of vector-Jacobian products. The cause of DENSERKS’s super-linear behavior is less clear.

\section{Conclusion\label{sec:conclusion}}
This article introduces a library designed for discrete adjoint sensitivity analysis of explicit Runge-Kutta methods, offering an efficient solution for computing sensitivities in systems of ODEs and related optimisation problems. The method implemented in this library leverages operator-overloading reverse automatic differentiation via the AADC framework of \cite{Matlogica}. This approach enables efficient evaluation of vector-Jacobian products, bypassing the limitations commonly encountered with black-box automatic differentiation techniques. By employing vectorization, the library can compute multiple vector-Jacobian products concurrently via SIMD vectorization, demonstrating an advantage over existing state-of-the-art approaches for adjoint sensitivity analysis of systems with many outputs.

Validation and experimentation were conducted to assess the correctness, usability, and efficiency of the proposed library compared to existing state-of-the-art software. The results demonstrate superior performance in the sensitivity analysis of the GLV system for large $N$, primarily due to the use of SIMD vectorization. Additionally, the checkpointing scheme in this implementation is less memory-intensive than that of FATODE, though only by a constant factor. For a more robust performance analysis, it would be advisable to simulate systems with the same $N$ but with different randomly generated parameters $A$ and $r$, and then take the average of the results.

Looking forward, there are several potential avenues for further development. The library could benefit from an option to switch between different modes of trajectory recording. Exploring other parallelization techniques, such as assigning a thread to compute sensitivities for each batch of output variables while using vectorization for concurrent vector-Jacobian product calculations, could also prove advantageous. To the best of the author's knowledge, this type of concurrency feature has not yet been implemented in existing sensitivity analysis software. Future enhancements could also include implementing an optimal check-pointing scheme similar to that of TSAdjoint, as well as expanding support for additional numerical integration methods, such as standard implicit and multi-step schemes. Adding functionality for second-order adjoint sensitivity analysis would further enhance the library's utility.

In conclusion, this work introduces a vectorization strategy for discrete adjoint sensitivity analysis, resulting in superior performance compared to some state-of-the-art software in particular scenarios involving few outputs and numerous parameters.

\section*{Code availability} 
The code is available from the GitHub repository \url{https://github.com/RuiMartins1996/VectorizedAdjoint}. The code has been tested under the \texttt{gcc} compiler.

\section*{Acknowledgements}
The first author is supported by FCT (``Funda\c{c}\~ao para a Ci\^encia e a Tecnologia'') under the PhD grant with reference 2022.13426.BD. The authors were supported by Portuguese funds through the CIDMA - Center for Research and Development in Mathematics and Applications, and the Portuguese Foundation for Science and Technology (``FCT--Funda\c{c}\~ao para a Ci\^encia e a Tecnologia''), within project UIDB/04106/2020 and UIDP/04106/2020. 

The first author wishes to express his gratitude to his PhD advisor, Dr. Uwe Khäler, for reviewing the manuscript and offering valuable constructive feedback. The authors also extend their thanks to Dr. Daniel Duffy for reviewing the manuscript and offering helpful criticism. The authors would also like to thank Dr. Hong Zhang for his availability and helpfulness in discussing the numerical results in this manuscript. His help and expertise in the field was very valuable in presenting good quality results. Finally, the authors would like to thank the people in their personal lives whose support and helpful suggestions contributed to the writing of this manuscript.

\section*{Declaration of generative AI and AI-assisted technologies in the writing process}
During the preparation of this work the authors used Chat-GPT to improve grammar, readability, and eloquence. The AI was not used to generate any scientific content. After using this tool/service, the author(s) reviewed and edited the content as needed and take full responsibility for the content of the publication.

\appendix

\section{Computing the partial derivatives\label{appendix::calculations}}
In this appendix we compute the partial derivatives shown in equation \ref{eq::partialderivatives1}, \ref{eq::partialderivatives2}, \ref{eq::partialderivatives3} and \ref{eq::partialderivatives4}, which are necessary for the full specification of the algorithm of discrete adjoint sensitivity analysis.

\subsection{Partial derivatives with respect to the intermediate variables}
It is useful to reacll the definition of the $m$th Runge-Kutta intermediate state $\mathbf{u}^{m,n}$:

\begin{equation}
\label{intermediate_state_x_m_n_definition}
u_i^{m,n}(\mathbf{w}^{0,n},\dots,\mathbf{w}^{m-1,n}) \\
= w_i^{0,n}+\Delta t^n \sum_{l=1}^{m-1} a_{ml} w_i^{l,n}, \qquad m = 1,\dots,s,
\end{equation}
and the definition of the intermediate states,

\begin{equation}
	\mathbf{w}^{m,n}=
	F\big(
	\mathbf{u}^{m,n},
	\alpha,
	t^{m,n}
	\big) \qquad  m = 1,\dots,s.
\end{equation}
Now, consider an infinitesimal displacement of the intermediate variable $w^{q,n}_j$, $w^{q,n}_j\to w^{q,n}_j+\Delta w^{q,n}_j$, and take the limit of $\frac{\Delta w^{m,n}_i}{\Delta w^{q,n}_j}$ when $\Delta w^{q,n}_j\to 0$. This limit corresponds to the sought after partial derivatives, and we proceed to analyse this limit for different values of the index $m$. The infinitesimal displacement on the $w^{q,n}_j$ intermediate variable induces a displacement on the $m$th Runge-Kutta intermediate state $\mathbf{u}^{m,n}$, namely:

\begin{align}
u^{m,n}_i(\mathbf{w}^{0,n},\dots,\mathbf{w}^{m-1,n})   &\to	u^{m,n}_i(\mathbf{w}^{0,n},\dots,w^{q,n}_j+\Delta w^{q,n}_j,\dots,\mathbf{w}^{m-1,n}) \nonumber
\\
&=u^{m,n}_i(\mathbf{w}^{0,n},\dots,\mathbf{w}^{m-1,n})+\Delta u^{m,n}_i.
\end{align}
In turn, the displacement $\Delta u^{m,n}_i$ is given by:

\begin{equation}
	\Delta u^{m,n}_i = \frac{\partial u^{m,n}_i(\mathbf{w}^{0,n},\dots,\mathbf{w}^{m-1,n})}{\partial w_j^{q,n}}\Delta w_j^{q,n}\text{ (no sum in $j$)},
\end{equation}
to first order in $\Delta w_j^{q,n}$, and the partial derivatives $\frac{\partial u^{m,n}_i}{\partial w_j^{q,n}}$ can be analytically computed from expression \ref{intermediate_state_x_m_n_definition},

\begin{equation}
	\frac{\partial u^{m,n}_i}{\partial w^{q,n}_j} =
	\begin{cases}
		\delta_{ij}
		& \text{if } q=0\\
		\\
		\Delta t^n a_{mq} \delta_{ij}
		& \text{if } q<m
	\end{cases}
\end{equation}

To compute the partial derivatives for $m=s+1$, we simply take into account the following expression,
\begin{equation}
     w^{s+1,n}_i =  w^{0,n}_i + \Delta t^n \sum_{l=1}^s b_l w_i^{l,n},
\end{equation}
and so we have:
\begin{equation}
	\frac{\partial w^{s+1,n}_i}{\partial w^{q,n}_j} =
	\begin{cases}
		\delta_{ij} & \text{if } q=0\\
		\delta_{ij} b_q \Delta t^n & \text{if } q=1,\dots,s\\
	\end{cases}.
\end{equation}
To compute the partial derivatives for any $m=1,\dots,s$, we need first to ascertain how the displacement on $w^{q,n}_j$ affects the intermediate variables $w_i^{m,n}$, namely:

\begin{equation}
w^{m,n}_i = F(\mathbf{u}^{m,n},t^{m,n}) \to
F_i(\mathbf{u}^{m,n}+\Delta\mathbf{u}^{m,n},t^{m,n}) =
w^{m,n}_i + \Delta w^{m,n}_i
\end{equation}
where, by using the chain rule, the nudge $\Delta w^{m,n}_i$ is given by:

\begin{align}
\Delta w^{m,n}_i &= F_i(\mathbf{u}^{m,n}+\Delta\mathbf{u}^{m,n},t^{m,n})-F_i(\mathbf{u}^{m,n},t^{m,n})\nonumber
 \\ 
&= J^x_{il}(\mathbf{u}^{m,n},t^{m,n})\frac{\partial u^{m,n}_l}{\partial w_j^{q,n}}\Delta w_j^{q,n} \qquad
 \text{(no sum in $j$)},
\end{align}
to first order in $\Delta w_j^{q,n}$. In the above expression,  $J^x(\mathbf{u}^{m,n},t^{m,n})$ denotes the $N\times N$ Jacobian matrix of $F$ taking $\alpha_k$ and $t^{m,n}$ as constants. The matrix is computed at the $m$th Runge-Kutta intermediate state $\mathbf{u}^{m,n}$, parameter vector $\alpha$ and intermediate time $t^{m,n}$:

\begin{equation}
	J^x(\mathbf{u},\alpha,t) =
	\begin{bmatrix}
		&\frac{\partial F_1(\mathbf{x},\alpha,t)}{\partial x_1} &\dots &\frac{\partial F_1(\mathbf{x},\alpha,t)}{\partial x_N}
		\\
		&\vdots  &\ddots &\vdots  
		\\
		&\frac{\partial F_N(\mathbf{x},\alpha,t)}{\partial x_1} &\dots &\frac{\partial F_N(\mathbf{x},\alpha,t)}{\partial x_N}
	\end{bmatrix},
\end{equation}
Finally, the partial derivatives $\frac{\partial w^{m,n}_i}{\partial w^{q,n}_j}$, for $m=1,\dots,s$, are given by the following limit:

\begin{equation}
	\frac{\partial w^{m,n}_i}{\partial w^{q,n}_j} = \lim_{\Delta w_j^{q,n}\to 0} \frac{\Delta w^{m,n}_i}{\Delta w_j^{q,n}} = J^x_{il}(\mathbf{u}^{m,n},\alpha,t^{m,n})\frac{\partial u^{m,n}_l}{\partial w_j^{q,n}}.
\end{equation}
In conclusion, we get the following expressions for the partial derivatives $\frac{\partial w^{m,n}_i}{\partial w^{q,n}_j}$ for $m=1,\dots,s$:

\begin{equation}
	\frac{\partial w^{m,n}_i}{\partial w^{q,n}_j} =
	\begin{cases}
		J^x_{ij}(\mathbf{x}^{m,n},\alpha,t^{m,n})
		& \text{if } q=0\\
		\\
		J^x_{ij}(\mathbf{x}^{m,n},\alpha,t^{m,n})
		\Delta t^n a_{mq}
		& \text{if } q<m
	\end{cases},
\end{equation}
To compute the partial derivatives for $m=0$, we take into account the following expression,
\begin{equation}
     w^{0,n}_i =  w^{s+1,n-1}_i
\end{equation}
and so we have:

\begin{equation}
	\frac{\partial w^{0,n}_i}{\partial w^{s+1,n-1}_j} = \delta_{ij}
\end{equation}

\subsection{Partial derivatives with respect to the parameters}
Consider an infinitesimal displacement of the $k$th parameter $\alpha_k \to \alpha_k+\Delta\alpha_k$, which induces a displacement on $w^{m,n}_i$ given by:

\begin{equation}
	\Delta w^{m,n}_i =F_i(\mathbf{x}^{m,n},\mathbf{\alpha}+\Delta\mathbf{\alpha},t^{m,n})-F_{\mathbf{\alpha},i}(\mathbf{x}^{m,n},\alpha,t^{m,n})= J^\alpha_{ik}\Delta \alpha_k \qquad \text{(no sum in $k$)},
\end{equation}
to first order in $\alpha_k$. Here, $J^\alpha$ corresponds to the $N\times P$ Jacobian matrix of $F$ taking all inputs except $\alpha$ as constants, computed at state $\mathbf{u}^{m,n}$, intermediate time $t^{m,n}$ and given parameters $\alpha$:

\begin{equation}
	J^\alpha(\mathbf{u},\alpha,t) = 
	\begin{bmatrix}
		&\frac{\partial F_{1}(\mathbf{u},\alpha,t)}{\partial \alpha_1} &\dots &\frac{\partial F_{1}(\mathbf{u},\alpha,t) }{\partial \alpha_M}
		\\
		&\vdots  &\ddots &\vdots  
		\\
		&\frac{\partial F_{N}(\mathbf{u},\alpha,t) }{\partial \alpha_1} &\dots &\frac{\partial F_{N}(\mathbf{u},\alpha,t) }{\partial \alpha_M}
	\end{bmatrix}.
\end{equation}
Finally, the partial derivatives $\frac{\partial w^{m,n}_i}{\partial \alpha_k}$ for $m=1,\dots,s$ are given by taking the following limit:

\begin{equation}
\frac{\partial w_i^{m,n}}{\partial \alpha_k} = \lim_{\Delta\alpha_k\to0} \frac{\Delta w_i^{m,n}}{\Delta\alpha_k}= J^\alpha_{ik}(\mathbf{u}^{m,n},\alpha,t^{m,n})
\end{equation}
and for $m=0$ and $m=s+1$ we have:

\begin{equation}
	\frac{\partial w_j^{0,n}}{\partial \alpha_k} = 0 \qquad\text{and}\qquad \frac{\partial w_j^{s+1,n}}{\partial \alpha_k} = 0
\end{equation}

 \bibliography{elsarticle}

\begin{thebibliography}{10}
\expandafter\ifx\csname url\endcsname\relax
  \def\url#1{\texttt{#1}}\fi
\expandafter\ifx\csname urlprefix\endcsname\relax\def\urlprefix{URL }\fi
\expandafter\ifx\csname href\endcsname\relax
  \def\href#1#2{#2} \def\path#1{#1}\fi

\bibitem{navon1998practical}
I.~M. Navon, Practical and theoretical aspects of adjoint parameter estimation and identifiability in meteorology and oceanography, Dynamics of Atmospheres and Oceans 27~(1-4) (1998) 55--79.

\bibitem{zhang2014parameter}
L.~Zhang, C.~Lyu, G.~Hinds, L.~Wang, W.~Luo, J.~Zheng, K.~Ma, Parameter sensitivity analysis of cylindrical lifepo4 battery performance using multi-physics modeling, Journal of The Electrochemical Society 161~(5) (2014) A762.

\bibitem{KPP2003}
V.~Damian, A.~Sandu, M.~Damian, F.~Potra, G.~R. Carmichael, The kinetic preprocessor kpp-a software environment for solving chemical kinetics, Computers \& Chemical Engineering 26~(11) (2002) 1567--1579.

\bibitem{sommer2017numerical}
A.~Sommer, Numerical methods for parameter estimation in dynamical systems with noise with applications in systems biology, Ph.D. thesis, Heidelberg University, Germany (2017).

\bibitem{sandu2005adjoint}
A.~Sandu, D.~N. Daescu, G.~R. Carmichael, T.~Chai, Adjoint sensitivity analysis of regional air quality models, Journal of Computational Physics 204~(1) (2005) 222--252.

\bibitem{ozyurt2005large}
D.~B. {\"O}zyurt, P.~I. Barton, Large-scale dynamic optimization using the directional second-order adjoint method, Industrial \& engineering chemistry research 44~(6) (2005) 1804--1811.

\bibitem{le2002second}
F.-X. Le~Dimet, I.~M. Navon, D.~N. Daescu, Second-order information in data assimilation, Monthly Weather Review 130~(3) (2002) 629--648.

\bibitem{griesse2003parametric}
R.~Griesse, A.~Walther, Parametric sensitivities for optimal control problems using automatic differentiation, Optimal Control Applications and Methods 24~(6) (2003) 297--314.

\bibitem{chen2018neural}
R.~T. Chen, Y.~Rubanova, J.~Bettencourt, D.~K. Duvenaud, Neural ordinary differential equations, Advances in neural information processing systems 31 (2018).

\bibitem{grathwohl2018ffjord}
W.~Grathwohl, R.~T. Chen, J.~Bettencourt, I.~Sutskever, D.~Duvenaud, Ffjord: Free-form continuous dynamics for scalable reversible generative models, arXiv preprint arXiv:1810.01367 (2018).

\bibitem{poli2019graph}
M.~Poli, S.~Massaroli, J.~Park, A.~Yamashita, H.~Asama, J.~Park, Graph neural ordinary differential equations, arXiv preprint arXiv:1911.07532 (2019).

\bibitem{burden2011numerical}
R.~L. Burden, Numerical analysis, Brooks/Cole Cengage Learning, 2011.

\bibitem{margossian2019review}
C.~C. Margossian, A review of automatic differentiation and its efficient implementation, Wiley interdisciplinary reviews: data mining and knowledge discovery 9~(4) (2019) e1305.

\bibitem{petsc2022}
H.~Zhang, E.~M. Constantinescu, B.~F. Smith, Petsc tsadjoint: a discrete adjoint ode solver for first-order and second-order sensitivity analysis, SIAM Journal on Scientific Computing 44~(1) (2022) C1--C24.

\bibitem{odessa1988}
J.~R. Leis, M.~A. Kramer, Algorithm 658: Odessa--an ordinary differential equation solver with explicit simultaneous sensitivity analysis, ACM Transactions on Mathematical Software (TOMS) 14~(1) (1988) 61--67.

\bibitem{cvodes2005}
R.~Serban, A.~C. Hindmarsh, Cvodes: the sensitivity-enabled ode solver in sundials, in: International Design Engineering Technical Conferences and Computers and Information in Engineering Conference, Vol. 47438, 2005, pp. 257--269.

\bibitem{daspk1999}
S.~Li, L.~Petzold, Design of new daspk for sensitivity analysis, UCSB Department of Computer Science Technical Report (1999).

\bibitem{juliaarticle}
Y.~Ma, V.~Dixit, M.~J. Innes, X.~Guo, C.~Rackauckas, A comparison of automatic differentiation and continuous sensitivity analysis for derivatives of differential equation solutions, in: 2021 IEEE High Performance Extreme Computing Conference (HPEC), IEEE, 2021, pp. 1--9.

\bibitem{DENSERKS2009}
M.~Alexe, A.~Sandu, Forward and adjoint sensitivity analysis with continuous explicit runge--kutta schemes, Applied Mathematics and Computation 208~(2) (2009) 328--346.

\bibitem{nixelotz2015higher}
J.~Lotz, U.~Naumann, R.~Hannemann-Tamas, T.~Ploch, A.~Mitsos, Higher-order discrete adjoint ode solver in c++ for dynamic optimization, Procedia Computer Science 51 (2015) 256--265.

\bibitem{fatode2014}
H.~Zhang, A.~Sandu, Fatode: A library for forward, adjoint, and tangent linear integration of odes, SIAM Journal on Scientific Computing 36~(5) (2014) C504--C523.

\bibitem{mitusch2019dolfin}
S.~K. Mitusch, S.~W. Funke, J.~S. Dokken, dolfin-adjoint 2018.1: automated adjoints for fenics and firedrake, Journal of Open Source Software 4~(38) (2019) 1292.

\bibitem{differentialequationsjl2017}
C.~Rackauckas, Q.~Nie, Differentialequations. jl--a performant and feature-rich ecosystem for solving differential equations in julia, Journal of open research software 5~(1) (2017).

\bibitem{rackauckas2020universal}
C.~Rackauckas, Y.~Ma, J.~Martensen, C.~Warner, K.~Zubov, R.~Supekar, D.~Skinner, A.~Ramadhan, Universal differential equations for scientific machine learning, arXiv preprint arXiv:2001.04385 (2020).

\bibitem{Matlogica}
D.~Goloubentsev, {Matlogica} website, developer of the library aadc, \url{https://matlogica.com/}, accessed: 2023-12-22 (2023).

\bibitem{lambert1991numerical}
J.~D. Lambert, et~al., Numerical methods for ordinary differential systems, Vol. 146, Wiley New York, 1991.

\bibitem{butcher2016numerical}
J.~C. Butcher, Numerical methods for ordinary differential equations, John Wiley \& Sons, 2016.

\bibitem{fehlberg1969low}
E.~Fehlberg, Low-order classical Runge-Kutta formulas with stepsize control and their application to some heat transfer problems, Vol. 315, National aeronautics and space administration, 1969.

\bibitem{odeint}
{odeint: a library for solving initial value problems}, \url{https://www.boost.org/doc/libs/1_82_0/libs/numeric/odeint/doc/html/index.html}, accessed: 2023-07-24 (2023).

\bibitem{zhuang2020ordinary}
J.~Zhuang, N.~Dvornek, X.~Li, J.~S. Duncan, Ordinary differential equations on graph networks, https://openreview.net/forum?id=SJg9z6VFDr (2019).

\bibitem{lotka2002contribution}
A.~J. Lotka, Contribution to the theory of periodic reactions, The Journal of Physical Chemistry 14~(3) (2002) 271--274.

\bibitem{volterra1926variazioni}
V.~Volterra, Variazioni e fluttuazioni del numero d'individui in specie animali conviventi, Societ{\`a} anonima tipografica" Leonardo da Vinci", 1926.

\bibitem{may1972will}
R.~M. May, Will a large complex system be stable?, Nature 238~(5364) (1972) 413--414.

\bibitem{Giering1998a}
R.~Giering, T.~Kaminski, {R}ecipes for {A}djoint {C}ode {C}onstruction, ACM Trans. On Math. Software 24~(4) (1998) 437--474.

\bibitem{hascoet2013tapenade}
L.~Hascoet, V.~Pascual, The tapenade automatic differentiation tool: principles, model, and specification, ACM Transactions on Mathematical Software (TOMS) 39~(3) (2013) 1--43.

\bibitem{petscautodiff}
J.~G. Wallwork, P.~Hovland, H.~Zhang, O.~Marin, Computing derivatives for petsc adjoint solvers using algorithmic differentiation, arXiv preprint arXiv:1909.02836 (2019).

\bibitem{adol1996}
A.~Griewank, D.~Juedes, J.~Utke, Algorithm 755: Adol-c: A package for the automatic differentiation of algorithms written in c/c++, ACM Transactions on Mathematical Software (TOMS) 22~(2) (1996) 131--167.

\end{thebibliography}
 \bibliographystyle{elsarticle}

\end{document}